\documentclass{amsart}

\usepackage{url,graphicx,amssymb,enumerate,stmaryrd, amsthm}
\usepackage[pagebackref,colorlinks,citecolor=blue,linkcolor=blue,urlcolor=blue,filecolor=blue]{hyperref}
\usepackage{microtype}
\usepackage[all]{xy}
\usepackage{dsfont}
\usepackage{aliascnt}
\usepackage[vcentermath]{youngtab}
\usepackage[dvipsnames, hyperref]{xcolor}
\usepackage{pinlabel}
\usepackage{etoolbox}
\usepackage{tikz} \usetikzlibrary{matrix,arrows}
\usepackage[hmargin=3cm,vmargin=3cm]{geometry}

\newcommand{\m}{\to}

\providecommand{\PBC}{\ensuremath\mathrm{PBC}}
\providecommand{\PB}{\ensuremath\mathrm{PB}}
\providecommand{\SPB}{\ensuremath\mathrm{SPB}}

\providecommand{\FI}{\ensuremath\mathsf{FI}}
\providecommand{\cI}{\ensuremath\mathcal{I}}
\providecommand{\VIC}{\ensuremath\mathsf{VIC}}

\providecommand{\SI}{\ensuremath\mathsf{SI}}
\providecommand{\AutF}{\ensuremath\mathsf{Aut}\, \ensuremath\mathsf{F}}
\providecommand{\C}{\ensuremath\mathcal{C}}
\providecommand{\G}{\ensuremath\mathcal{G}}

\providecommand{\hooklongrightarrow}{\lhook\joinrel\longrightarrow}
\providecommand{\twoheadlongrightarrow}{\relbar\joinrel\twoheadrightarrow}
\providecommand{\inject}{\hooklongrightarrow}
\providecommand{\surject}{\twoheadlongrightarrow}

\providecommand{\N}{\ensuremath\mathbb N_0}
\providecommand{\Z}{\ensuremath\mathbb Z}
\providecommand{\Q}{\ensuremath\mathbb Q}

\DeclareMathOperator{\im}{im}
\DeclareMathOperator{\coker}{coker}
\DeclareMathOperator{\coeq}{coeq}
\DeclareMathOperator{\Hom}{Hom}

\providecommand{\id}{\ensuremath\mathrm{id}}

\DeclareMathOperator{\Aut}{Aut}

\DeclareMathOperator{\GL}{GL}

\DeclareMathOperator{\Sp}{Sp}

\DeclareMathOperator{\IA}{IA}
\DeclareMathOperator{\Mod}{Mod}

\DeclareMathOperator{\ord}{ord}
\DeclareMathOperator{\Lk}{Lk}
\DeclareMathOperator{\Ind}{Ind}
\DeclareMathOperator{\Res}{Res}

\DeclareMathOperator{\rk}{rk}

\DeclareMathOperator{\ab}{ab}

\DeclareMathOperator{\spn}{span}
\DeclareMathOperator{\sat}{sat}
\newcommand{\U}{\mathfrak U}

\definecolor{grey}{gray}{.5}

\providecommand{\con}[1]{\textbf{\textup{#1}}}

\numberwithin{thmcounter}{section}
\newaliascnt{thmauto}{thmcounter}

\newaliascnt{Defauto}{thmcounter}

\newaliascnt{exauto}{thmcounter}

\newaliascnt{lemauto}{thmcounter}

\newaliascnt{propauto}{thmcounter}

\newaliascnt{corauto}{thmcounter}

\newaliascnt{remauto}{thmcounter}

\newaliascnt{convauto}{thmcounter}

\newtheorem{atheorem}{Theorem}

\newtheorem*{ThmA'}{Theorem A'}
\newtheorem*{ThmB'}{Theorem B'}
\newtheorem*{ThmC'}{Theorem C'}

\newtheorem{theorem}[thmauto]{Theorem}
\newtheorem{lemma}[lemauto]{Lemma}
\newtheorem{proposition}[propauto]{Proposition}
\newtheorem{corollary}[corauto]{Corollary}
\theoremstyle{definition}
\newtheorem{definition}[Defauto]{Definition}
\newtheorem{example}[exauto]{Example}
\newtheorem{remark}[remauto]{Remark}
\newtheorem{convention}[convauto]{Convention}

\title[Central stability for congruence subgroups and Torelli groups]{Central stability for the homology of congruence subgroups and the second homology of Torelli groups}
\author{Jeremy Miller}\thanks{Jeremy Miller was supported in part by NSF grant DMS-1709726}
\address{Department of Mathematics, Purdue University, USA}
\email{jeremykmiller@purdue.edu}
\author{Peter Patzt}
\address{Department of Mathematics, Purdue University, USA}
\email{ppatzt@purdue.edu}
\author{Jennifer C. H. Wilson}
\address{Department of Mathematics, University of Michigan, USA}
\email{jchw@umich.edu}
\date{September 2016}

 \date{\today}

\begin{document}

\begin{abstract}
We prove a representation stability result for the second homology groups of Torelli subgroups of mapping class groups and automorphism groups of free groups. This  strengthens the results of Boldsen--Hauge Dollerup and Day--Putman. We also prove a new representation stability result for the homology of certain congruence subgroups, partially improving upon the work of Putman--Sam. These results follow from a general theorem on syzygies of certain modules with finite polynomial degree.

\end{abstract}
\maketitle

\tableofcontents

%%%%%%%%%%%%% SECTION %%%%%%%%%%%%%%%%%%%%

\section{Introduction}

The purpose of this paper is to prove cases of a variant of two conjectures from Church and Farb \cite[variations on Conjectures 6.3 and 6.5]{CF} as well as to answer a question posed by Putman \cite[fifth Remark]{Pu}. Specifically, we prove a central stability result for the second homology groups of the Torelli subgroups of mapping class groups and automorphism groups of free groups, and a central stability result for the homology groups in all degrees of certain congruence subgroups of general linear groups. We begin by recalling the definition of central stability.  

%These results improve some of the results of Boldsen and Hauge Dollerup \cite[Theorem 1.0.1]{BHD}, Day and Putman \cite[Theorem B]{DP}, and Putman and Sam \cite{PutmanSamLinearGroups}.  

%The main technical lemma of this paper is a result showing higher central stability for modules over certain categories with finite polynomial degree, such as modules over the categories appearing in work of Putman and Sam \cite{PutmanSamLinearGroups}. The proof of this lemma involves verifying a general acyclicity criterion  developed by the second author \cite{Pa2} for complexes associated to these modules. From this result we deduce that if a module satisfies a polynomial condition -- a condition often appearing in twisted homological stability theorems -- then this module exhibits a form of representation stability.

\subsection{Central stability}

 Let \[G_0 \longrightarrow G_1 \longrightarrow G_2 \longrightarrow \cdots\]
be a sequence of groups and group homomorphisms, and let \[A_0 \xrightarrow{\phi_0} A_1 \xrightarrow{\phi_1} A_2 \xrightarrow{\phi_2} \cdots\]
be a sequence with $A_n$ a $G_n$--representations and $\phi_n$ a $G_n$--equivariant map. In this paper, the groups $G_n$ will be symmetric groups, general linear groups, subgroups of general linear groups with restricted determinant, or symplectic groups. Church and Farb \cite{CF} described a condition they called \emph{representation stability} for these types of sequences $\{A_n\}$. In its original formulation, representation stability was only defined for semisimple representations $A_n$ and for families of groups admitting natural identifications between irreducible $G_n$--representations for different $n$. A much more general and formal notion of representation stability called \emph{central stability} was introduced by Putman \cite{Pu}.

 For purposes of exposition, we will now specialize to the case that the groups $G_n$ are either the family of general linear groups $\GL_n(\Z)$ or the family of symmetric groups $\mathfrak S_n$. Let $\sigma_n \in G_n$ be the transposition $(n-1\;n)$ in the case $G_n=\mathfrak S_n$ and the associated permutation matrix in the case $G_n=\GL_n(\Z)$.  We assume that $\sigma_n$ fixes %the image of 
the map $(\phi_{n-1} \circ \phi_{n-2})\colon A_{n-2} \to A_{n}$. There are two natural maps
  \[\Ind_{G_{n-2}}^{G_n} A_{n-2} \rightrightarrows  \Ind_{G_{n-1}}^{G_n} A_{n-1}, \] 
the first induced by $\phi_{n-2}$ and the second by postcomposing this induced map by  $\sigma_{n}$.  We say that the sequence $\{A_n\}$ has \emph{central stability degree} $\le d$ if for all $n> d$ the map 
\[ \coeq \left(\Ind_{G_{n-2}}^{G_n} A_{n-2} \rightrightarrows  \Ind_{G_{n-1}}^{G_n} A_{n-1} \right) \m A_n \]
induced by $\phi_{n-1}$ is an isomorphism. %such that $\phi_{n-1}$ coincides with the natural map from $A_{n-1}$ to this coequalizer.
We say that $\{A_n\}$ is \emph{centrally stable} if it has finite central stability degree. If $\{A_n\}$ has central stability degree $\le d$, then the entire sequence is determined by the finite sequence $A_0 \to A_1 \to \cdots \to A_{d}$. Analogous definitions exist for symplectic groups and subgroups of the general linear groups with restricted determinant, which we review in \autoref{sec3}. The main result of this paper is to prove that certain homology groups of Torelli groups and congruence subgroups are centrally stable. 

A cautionary remark for experts: there are different definitions of \emph{central stability} in the literature that are not equivalent in all situations; see \autoref{RemarkCSWarning}. 

\filbreak
\subsection{Central stability for Torelli groups and congruence subgroups}
\subsubsection*{Automorphisms of free groups} Let $F_n$ denote the free group on $n$ letters. A group automorphism $f\colon F_n \to F_n$ induces a linear map $f^{\ab}\colon \Z^n \to \Z^n$ on the abelianization of $F_n$. This construction defines a surjective map \[\Aut(F_n) \to \GL_n(\Z).\] We write $\IA_n$ to denote the kernel of this map, which is often called the \emph{Torelli subgroup} of $\Aut(F_n)$.  If we fix an inclusion of $\Aut(F_n)$ into $\Aut(F_{n+1})$ as the stabilizer of a free generator, the Torelli group $\IA_n$ maps into $\IA_{n+1}$. The induced maps 
$H_i(\IA_n) \to H_i(\IA_{n+1})$ are equivariant with respect to the induced actions of the groups $\GL_n(\Z)$, and so we may ask whether  the sequences of homology groups $\{H_i(\IA_n)\}$ exhibit central stability as $\GL_n(\Z)$--representations.

The first homology group of $\IA_n$  was computed independently in work of Andreadakis \cite[see Section 6]{Andr} for $n=2,3$, Cohen--Pakianathan (unpublished), Farb (unpublished) and Kawazumi \cite[Theorem 6.1]{Kawazumi}. Very little is known about the higher homology groups of $\IA_n$. Bestvina--Bux--Margalit \cite[Main Theorem, part (2)]{BBM} proved that $H_2(\IA_3)$ is not finitely generated as an abelian group which implies that $\IA_3$ is not finitely presented, a result originally due to Krsti\'c--McCool \cite[Theorem 1]{KMc}. However, for $n>3$, it is unknown if $H_2(\IA_n)$ is finitely generated or if $\IA_n$ is  finitely presented.

%Even in degree two, however, the homology groups remain fairly mysterious.

Church and Farb  \cite[Conjecture 6.3]{CF} conjectured ``mixed representation stability'' for a certain summand of the rational homology of $\IA_n$. Given Putman's subsequent work on central stability, it is natural to modify Church and Farb's Conjecture 6.3 to also ask if the integral homology groups of $\IA_n$ exhibit central stability as representations of the general linear groups. We prove this result for the second homology groups.

\begin{atheorem} \label{IAn}
The sequence $H_2(\IA_n)$ has central stability degree $\le 38$ as $\GL_n(\Z)$--representations. 
%\todo{Fix bounds.}
\end{atheorem}

This theorem partially improves upon a result of Day and Putman \cite[Theorem B]{DP}, which established surjectivity of the maps \[\coeq \left(\Ind_{\GL_{n-2}(\Z)}^{\GL_{n}(\Z)} H_2(\IA_{n-2}) \rightrightarrows  \Ind_{\GL_{n-1}(\Z)}^{\GL_{n}(\Z)} H_2(\IA_{n-1}) \right) \to H_2(\IA_{n})\]   for $n>6$. Our techniques only show these maps surject for $n>18$ (\hyperlink{thmA'}{Theorem A'}), but additionally prove the maps are injective for  $n>38$.  We note that our method is substantially different from that of Day and Putman: our proof centers on properties of general linear groups, while their proof focused on properties of automorphism groups of free groups. Our proof strategy applies to Torelli subgroups of mapping class groups with little modification. In contrast, Day and Putman \cite[Remark 1.3]{DP} noted that their techniques for Torelli subgroups of $\Aut(F_n)$ do not easily generalize to Torelli subgroups of mapping class groups, and that the techniques used by Boldsen and Hauge Dollerup \cite{BHD} for Torelli subgroups of mapping class groups do not easily generalize to prove results about $\IA_n$.

{   %%%%%%%%%%%%%%%%%%%%%%%% BEGIN PURPLE
Our proof uses a spectral sequence due to Putman--Sam \cite{PutmanSamLinearGroups}, %(\autoref{Def:stability SES}), 
and relies on known computations of $H_1(\IA_n)$ (Kawazumi \cite[Theorem 6.1]{Kawazumi}) and  a high-connectivity result of Hatcher--Vogtmann \cite[Proposition 6.4]{HatcherVogtmannCerf}. 
} %%%%%%%%%%%%%%%%%%%%%%%% END PURPLE

\subsubsection*{Mapping class groups}  Let $\Sigma_{g,1}$ denote an oriented genus-$g$ surface with one boundary component, and let $\Mod_g$ denote the group of connected components of the group of orientation-preserving diffeomorphisms of $\Sigma_{g,1}$ that fix the boundary pointwise. The induced action of these diffeomorphisms on the first homology group of $\Sigma_{g,1}$ preserves the intersection form, which happens to be symplectic. This action therefore induces a map \[\Mod_g \to \Sp_{2g}(\Z),\] where $\Sp_{2g}(\Z)$ denotes the group of linear automorphisms of $\Z^{2g}$ that preserve the standard symplectic form. This map is well known to be surjective. The kernel $\mathcal I_g$ of this map is called the \emph{Torelli subgroup} of $\Mod_g$. Just as with the groups $\IA_n$, very little is known about the homology groups $H_i(\mathcal I_g)$ for $i>1$. Church and Farb conjectured a form of representation stability for $H_i(\mathcal I_g)$ \cite[Conjecture 6.1]{CF}, motivated by their earlier work \cite{CFAbelianCycles} constructing nontrivial classes in $H_i(\mathcal I_g ;\Q)$. In this paper we prove a central stability result for the homology groups $H_2(\mathcal I_g)$ as symplectic group representations.

\begin{atheorem} \label{Ig}
The sequence $H_2(\mathcal I_g)$ has central stability degree $\le 69$ as $\Sp_{2g}(\Z)$--representations.
\end{atheorem}

This partially improves upon a result of Boldsen and Hauge Dollerup \cite[Theorem 1.0.1]{BHD} proving that \[\coeq \left(\Ind_{\Sp_{2g-4}(\Z)}^{\Sp_{2g}(\Z)} H_2(\mathcal I_{g-2};\Q) \rightrightarrows  \Ind_{\Sp_{2g-2}(\Z)}^{\Sp_{2g}(\Z)} H_2(\mathcal I_{g-1};\Q) \right) \to H_2(\mathcal I_{g};\Q)\] is surjective for $g>6$; our techniques only show the map is surjective for $g>33$ (\hyperlink{thmB'}{Theorem B'}) but also establish injectivity for $g>69$, and our result holds with integral coefficients. 

{  %%%%%%%%%%%%%%%%%%%%%%%% BEGIN PURPLE
Our proof again uses the spectral sequence from Putman--Sam \cite{PutmanSamLinearGroups}, % (see \autoref{Def:stability SES}), 
 Johnson's computations of $H_1(\mathcal{I}_g)$ (\cite[Theorem 3]{JohnsonAbelianization}) and  a high--connectivity result of Hatcher--Vogtmann  \cite[Main Theorem]{HatcherVogtmannTethers}.  
} %%%%%%%%%%%%%%%%%%%%%%%% END PURPLE

Prior work had established representation stability results for certain subrepresentations of the rational cohomology of Torelli groups. Hain \cite[Remark 14.2 and Proposition 5.1]{Hain} computed the images of the cup product maps
          $$H^1(\mathcal I_g; \Q) \otimes H^1(\mathcal I_g;\Q) \longrightarrow H^2(\mathcal I_g;\Q)$$
as $\Sp_{2g}(\Z)$-representations, and his results show in particular that the multiplicities stabilize for $g>5$. Church--Ellenberg--Farb \cite[proof of Theorem 7.2.2]{CEF} proved a representation stability result for the subalgebras of $H^*(\mathcal I_g; \Q)$ generated by $H^1(\mathcal I_g; \Q)$, viewed as $S_n$-representations. They also proved the analogous result for $\IA_n$ \cite[proof of Theorem 7.2.3]{CEF}. Kupers--Randal-Williams \cite[Theorem 8.1]{KRW} computed the algebraic part of $H^2(\mathcal I_g;\Q)$ for $g \geq 4$.

Shortly after this paper was released, Kassabov--Putman \cite{KPu} proved that $H_2(\mathcal I_g;\Z)$ is finitely generated as an $\Sp_{2g}
(\Z)$--module for $g \geq 3$ \cite[Theorem A]{KPu}, and that $\mathcal I_g$ has a \emph{finite $\Mod_g$-equivariant presentation}  for $g \geq 3$ \cite[Theorem C]{KPu}. 

\subsubsection*{Congruence subgroups of $\GL_n(R)$} We also investigate representation stability for congruence subgroups of general linear groups. Homological and representation stability properties of congruence subgroups have had applications in homotopy theory -- for example, in work of Charney \cite{Ch1} on excision for algebraic $K$-theory -- and applications in number theory; see Calegari--Emerton \cite{CaEm1}. Let $I$ be a two-sided ideal in a ring $R$, and let $\GL_n(R,I)$ denote the kernel of the ``reduction modulo $I$'' map
\[\GL_n(R) \to \GL_n(R/I).\] We call $\GL_n(R,I)$ the \emph{level $I$ congruence subgroup} of $\GL_n(R)$. 

%We note that, up to isomorphism, the group $\GL_n(R,I)$ only depends on the ideal $I$, viewed as a non-unital ring, and not the ambient ring $R$. %  Given a non-unital ring $I$ we can, in fact,  construct $\GL_n(R,I)$ abstractly, without reference to $R$: the group consists of $n \times n$ matrices whose off-diagonal entries are elements of $I$. The diagonal entries are formal objects of the form $(1+i)$ with $i \in I$ and $1$ a ``unit". Matrix multiplication is determined by the arithmetic laws in $I$.

 For $R$ commutative and $H$ a subgroup of the group of units $R^\times$ of $R$, let $\GL_n^H(R)$ denote the subgroup of matrices with determinant in $H$. Let $\U$ denote the image of $R^\times$ in $R/I$. In many cases, $\GL_n^\U(R/I)$ is the image of $\GL_n(R)$ in $\GL_n(R/I)$. For example, this is the case when $R$ is the ring of integers in a number field by strong approximation (see e.g \cite[Chapter 7]{PR}) or when $R/I$ is Euclidean since then $\GL_n^\U(R/I)$ will be generated by elementary matrices. In these cases, the homology groups $H_i(\GL_n(R,I))$ will have natural linear $\GL_n^{\mathfrak U}(R/I)$--actions and there are equivariant maps $H_i(\GL_n(R,I)) \to H_i(\GL_{n+1}(R,I))$ induced by the inclusions $\GL_n(R,I) \hookrightarrow \GL_{n+1}(R,I)$ .

%For a ring $I$ without unit, we denote by $I_+$ its \emph{unitalization} which is $\Z \oplus I$ as an abelian group and its multiplication is defined by
%\[ (m,x) \cdot (n,y) = (mn, nx+my+xy)\]
%for $(m,x),(n,y) \in I_+$. Then $I$ is a two-sided ideal in $I_+$.  Thus $\GL_n(R,I)$ is the level $I$ congruence subgroup of $\GL_n(I_+)$ and $\GL_n(\Z)$ acts linearly on $H_i(\GL_n(R,I))$. For every ring $R$ that contains $I$ as a two-sided ideal, there is a unique ring map $I_+ \to R$ that sends  is the identity on $I$.

%This was first done by Charney \cite{Ch1} and has recently received more attention in the context of representation stability.

%The group $\GL_n^\U(R/I)$ acts on $\GL_n(R,I)$.

When $I \subsetneq R$, the symmetric group $\mathfrak S_n$ is naturally a subgroup of $\GL_n^\U(R/I)$, so we may view the the groups $H_i(\GL_n(R,I))$ as representations of the symmetric group. When $I=R$, we endow $H_i(\GL_n(R,I))$ with the trivial symmetric group action. Putman  \cite[Theorem B]{Pu} proved that when $R$ has finite stable rank (see \autoref{DefnStableRank} and Bass  \cite[Section 4, ``Definition"]{BassKTheory}),  these homology groups have central stability as $\mathfrak S_n$--representations. Putman's result gave explicit stable ranges for the homology groups $H_i(\GL_n(R,I); K)$ over certain fields $K$, and later Church--Ellenberg--Farb--Nagpal \cite[Theorem D]{CEFN} proved a central stability result for the homology of certain congruence subgroups with coefficients in a general Noetherian ring $K$ but without explicit bounds. Church and Ellenberg \cite[Theorem D']{CE} generalized both theorems with a result for integral homology with explicit stable ranges, and Church--Miller--Nagpal--Reinhold \cite[Application B]{CMNR} and Gan--Li \cite[Theorem 11]{GanLiCongruenceSubgroups} further improved these ranges. 

%Although $\GL_n(R,I)$ is independent of $R$ as a representation of $\mathfrak S_n$ all results in \cite{Ch1, Pu, CEFN,CE} depend on the stable rank of $R$. Ideally these results would only depend on $I$. We have no solution for this. 

Putman \cite[fifth Remark]{Pu} comments that it would be ideal to understand stability properties of $\GL_n(R,I)$ as $\GL^{\mathfrak U}_n(R/I)$--representations instead of just $\mathfrak S_n$--representations. Using Gan--Li \cite[Theorem 11]{GanLiCongruenceSubgroups} (stated in our \autoref{CEcong}), we provide the following partial solution to this problem. 

%As we have pointed out before, though, it is no less natural to understand these homology groups as $\GL_n(\Z)$--representations.

\begin{atheorem}\label{GL(R,I)}
Let $I$ be a two-sided ideal in a ring $R$ of stable rank $r$ and let $t$ be the minimal stable rank of all rings containing $I$ as a two-sided ideal. Let $\U$ denote the image of $R^\times$ in $R/I$. If $R/I$ is a PID of stable rank $s$, and the natural map $ \GL_n(R) \m \GL^{\U}_n(R/I)  $ is surjective,  then the central stability degree of the sequence 
 $H_i(\GL_n(R,I))$ is \[\begin{array}{lr}
\le s+1 & \text{ for } i=0 \\
 \le  \min( 24i +12t +s +4, \, \max(5+r,5+s) ) & \text{ for } i=1 \\
 \le 24i +12t +s +4 & \text{ for } i>1
\end{array} \] as $\Z[\GL_n^\U(R/I)]$--modules. 

 %$H_0(\GL_n(R,I))$ has central stability degree \[\le s,\] $H_1(\GL_n(R,I))$ has central stability degree \[\le \max(4+t,6+s),\] and $H_i(\GL_n(R,I))$ has central stability degree  \[\le 3(2^{i-1})(2t+7)-8+s \] for $i \ge 2$ as $\GL_n(R/I)$--representations. 

%\[\le \begin{cases}
%s & \text{if } i=0 \\
%\max(4+t,6+s) & \text{if } i=1 \\
%3(2^{i-1})(2t+7)-8+s  & \text{if } i \ge 2 
%\end{cases}\]

\end{atheorem}

%{\bf \color{Purple} %%%%%%%%%%%%%%%%%%%%%%%%%%%%%%%%%%% BEGIN PURPLE
%\begin{atheorem}\label{GL(R,I)} NEED TO REPHRASE!!!!!!! POSSIBLY CAN LET $R$ BE NONCOMMUTATIVE BUT THEN NEED TO MAKE SURE THAT $\GL_n^\U(R/I)$ IS THE RIGHT GROUP. POSSIBLY CANNOT LET $t$ BE MIN BECAUSE NEED BOUNDS ON $R$ and $R/I$ TO BE COMPATIBLE. 
%Let $I$ be a proper nonzero two-sided ideal of a commutative ring $R$ and let $t$ be the minimal stable rank of all commutative rings containing $I$ as a proper ideal. If $R/I$ is a PID of stable rank $s$, then the central stability degree of the sequence 
% $H_i(\GL_n(R,I))$ is \[\begin{array}{lr}
%\le s+1 & \text{ for } i=0 \\
% \le \max(5+t,5+s)  & \text{ for } i=1 \\
% \le 24i +12t +s +4 & \text{ for } i>1
%\end{array} \] as $\Z[\GL_n^\U(R/I)]$--modules. 

%\end{atheorem}
%} %%%%%%%%%%%%%%%%%%%%%%%%%%%%%%%%%%% END PURPLE
 
 %%%%%%%%%%%%%%%%%%%%%%%%%%%%%%%%%%% BEGIN PURPLE

 %%%%%%%%%%%%%%%%%%%%%%%%%%%%%%%%%%% END PURPLE

{\color{Red} %%%%%%%%%%%%%%%%%%%%%%%%%%%%%%%%%%% BEGIN RED
%\begin{remark}
 %Note that the quotient $R/I$ need not be a PID in order to define an action on $H_*(\GL_n(R,I))$ by a group $\GL_n(Q)$ with $Q$ a PID. For example, given an non-unital ring $I$, we can realize $I$ as a two-sided ideal in its \emph{unitilization $I_+$}, its image under the left adjoint to the forgetful functor from the category of (unital) rings to the category of non-unital rings. Then $I_+/I \cong \Z$, and $\GL_n(I_+,I) \cong \GL_n(R,I)$. Hence in particular $H_*(\GL_n(R,I))$ is always a $\GL_n(\Z)$--representation. 
%\end{remark}
} %%%%%%%%%%%%%%%%%%%%%%%%%%%%%%%%%%% END RED

%\todo[Peter]{We should get a similar statement for $\Sp$.}

In particular, \autoref{GL(R,I)} applies in the case that $R$ is a ring of stable rank $t=r$, and $R/I$ is a Euclidean domain. It follows from the work of Putman and Sam \cite{PutmanSamLinearGroups} that $H_i(\GL_n(R,I))$ has finite central stability degree when $R$ is the ring of integers in an algebraic number field. However, their techniques do not give an explicit central stability range or address congruence subgroups with infinite quotients $R/I$. 

\begin{remark}\label{remarkUnitalization}
Observe that the group $\GL(R,I)$ only depends on $I$ as a non-unital ring. The only assumption we put on the non-unital ring $I$ in \autoref{GL(R,I)} is that $I$ is a two-sided ideal in a ring with finite Bass stable rank.   We can always find a ring $R$ containing $I$ as a two sided ideal with $ \GL_n(R) \m \GL^{\U}_n(R/I)  $  surjective. In particular, $R$ can be taken to be the \emph{unitalization} $I_+$, $R/I$ can be taken to be $\Z$, and $\U$ can be taken to be $\{\pm 1\}$. The unitalization of $I$ is defined as the image of $I$ under the left adjoint of the forgetful functor from the category of unital rings to the category of possibly non-unital rings. As an abelian group, $I_+ \cong \Z \oplus I$. The unitalization is an important construction when studying excision for algebraic $K$-theory (see e.g. Suslin--Wodzicki \cite{SW}). In particular, for any non-unital ring $I$, $H_i(\GL_n(I_+,I))$ is a $\GL_n(\Z)$--representation. 

Because $\GL(R,I) \cong \GL_n(I_+,I)$ for any ring $R$ and ideal $I$, and because $\Z$ has stable rank $2$, we obtain the following corollary of \autoref{GL(R,I)}. 

\begin{corollary}

Let $I$ be a two-sided ideal in a ring $R$ and let $t$ be the minimal stable rank of all rings containing $I$ as a two-sided ideal.   Then, as a sequence of $\Z[\GL_n(\Z)]$--modules,  the sequence $H_i(\GL_n(R,I))$ has central stability degree \[\begin{array}{lr}

\le 3 & \text{ for } i=0 \\

 \le  \min( 24i +12t +6, \, 7 ) & \text{ for } i=1 \\

 \le 24i +12t +6 & \text{ for } i>1. 

\end{array} \] 

\end{corollary}

\end{remark}

\begin{remark}

Currently, the representation stability literature has a number of effective tools for establishing stability results, including central stability, for sequences of symmetric group representations (see for example Church--Ellenberg \cite{CE}) but often these tools cannot be applied to sequences of general linear (or symplectic) group representations.  In \autoref{Cor:FIvsVICH}, we prove that if a sequence of general linear representations happens to be centrally stable with respect to the underlying actions of the symmetric groups, then it is also centrally stable with respect to the action of the general linear groups (though the converse does not hold in general). \autoref{Cor:FIvsVICH} then allows us to use known or elementary results about central stability for certain sequences of symmetric group representations to prove our central stability results for general linear group representations, in particular \autoref{IAn} and \autoref{GL(R,I)}.

We prove \autoref{Cor:FIvsVICH} by comparing central stability with the notion of \emph{polynomial degree} (see \autoref{defpoly}). This polynomial condition and its variants have a long history in the literature; see e.g. \cite{EilenbergMacLaneII, Pirashvili-PolynomialFunctors, Pirashvili-FiniteFields, BauesPirashvili, Pirashvili-DoldKan, BauesDreckmannFranjouPirashvili,  HartlVespa-Quadratic, DjamentUnitaires,  DjamentVespa-tordus, VespaThesis,  DjamentVespa-hermitiens, DjamentVespa-WeakPolynomial, DjamentPirashviliVespa, Djament-HomologieStable, Djament-finitude,RWW}. The definition used in this paper agrees with that of an arXiv version of Randal-Williams--Wahl \cite[Definition 4.10]{RWWarxiv}. A convenient fact is that the polynomial degree of a sequence does not depend on the choice of  automorphism group, which allows us to leverage results about symmetric group actions to prove results about general linear group actions. We use a similar strategy in the symplectic case as well in order to prove \autoref{Ig}.

We note that we do not prove that the second homology groups of these Torelli groups satisfy a polynomial condition. Instead, we use that the zeroth and first homology groups do statisfy a polynomial condition to prove that the second homology groups are centrally stable. This contrasts with the case of congruence subgroups where we use that all of the homology groups satisfy a polynomial condition.

%There is a vast literature on central stability for $\FI$-modules. By viewing the symmetric group as the group of permutation matrices, we can view a $\VIC(R)$--module $A$ as an $\FI$--module. The polynomial degree of $A$ as a $\VIC(R)$--module agrees with its polynomial degree as an $\FI$--module. Using previously known results or elementary calculations, it is often easy to prove central stability for $A$ as an $\FI$-module. This lets us bound the polynomial degree of $A$ which in turn lets us prove (higher) central stability of $A$ as a $\VIC(R)$--module. This strategy plays a fundamental role in all of our applications.

\end{remark}

\subsection{Outline}

In \autoref{sechighcon}, we prove that some relevant semisimplicial sets are highly connected. We use these connectivity results in \autoref{sec3} to prove that certain modules with finite polynomial degree exhibit central stability.  In \autoref{secap}, we apply our results on polynomial degree to prove \autoref{IAn}, \autoref{Ig}, and \autoref{GL(R,I)}.

\subsection{Acknowledgments} 

This collaboration started as a result of the 2016 AIM Workshop on Representation Stability. We would like to thank AIM as well as the organizers of the workshop, Andrew Putman, Steven Sam, Andrew Snowden, and David Speyer. Additionally, we thank Benson Farb, Thomas Church, Dick Hain, Alexander Kupers, Andrew Putman, Holger Reich, Steven Sam, and Graham White for helpful conversations and feedback. We are grateful to the referees of this paper for their careful reading and detailed feedback. 

%\todo{Re-institute acknowledgments}

%%%%%%%%%%%%% SECTION %%%%%%%%%%%%%%%%%%%%

\section{High connectivity results}    \label{sechighcon}

As is common in stability arguments, our proofs will involve establishing high connectivity for certain spaces with actions of our families of groups. We start out in  \autoref{SecSimplicialReview} with a review of simplicial complexes and semisimplicial sets, and an overview of techniques to prove that their geometric realizations are highly connected. Then in  \autoref{SecAlgebraReview} we review elementary properties of free modules and symplectic structures over PIDs. In \autoref{SecPB} we prove some new high connectivity results for spaces relevant to stability for general linear group representations. In \autoref{SecSPB}, we review results from the literature concerning high connectivity of spaces relevant to stability for symplectic group representations. 

\subsection{Review of simplicial techniques} \label{SecSimplicialReview} 

 %%%%%%%%%%%%% START PURPLE
Recall that the data of a semisimplicial set is the same as the data of a simplicial set without degeneracy maps; see Weibel \cite[Definition 8.1.9]{Weibel}. %Randal-Williams--Wahl \cite[Section 2.1]{RWW} gives a comparison of semisimplicial sets versus simplicial complexes.
 %%%%%%%%%%%%%%%%%%%%%%% END PURPLE
 We will also consider simplicial complexes in this paper. Simplicial complexes differ from semisimplicial sets in several ways. For example, each simplex of a semisimplicial set comes equipped with an order on its set of faces, while the faces of a simplicial complex are not ordered. Further, a collection of vertices can be the set of vertices of at most one simplex of a simplicial complex, but there are no such restrictions for semisimplicial sets. See Randal-Williams--Wahl \cite[Section 2.1]{RWW} for a discussion of the differences between semisimplicial sets and simplicial complexes. 

We say that a semisimplicial set or simplicial complex is $n$--connected if its geometric realization is $n$--connected. If $\sigma$ is a simplex of a simplicial complex $X_\circ$, we write $\Lk^X(\sigma)$ to denote the link of $\sigma$ in $X$, or simply $\Lk(\sigma)$ when the ambient complex is clear from context. We now recall the definition of weakly Cohen--Macaulay simplicial complexes. 

\begin{definition}
A simplicial complex $X_\circ$ is called \emph{weakly Cohen--Macaulay} (abbreviated wCM) of dimension $n$ if it satisfies the following two conditions.
\begin{itemize}
  \item  $X_\circ$ is $(n-1)$--connected.
  \item If $\sigma$ is a $p$--simplex of $X_\circ$, then $\Lk^X(\sigma)$ is $(n-2-p)$--connected.
\end{itemize} 
\end{definition}

\begin{definition}
If $X_\circ$ is a simplicial complex, let $X_\bullet = X^{\ord}_\circ$ be the associated semisimplicial set formed by taking all simplices of $X_\circ$ with all choices of orderings on their vertices.  \end{definition}

\begin{convention} We adopt the following convention on subscript notation: we always write  $X_\bullet$  for a semisimplicial set and $X_\circ$  for a simplicial complex.  If we denote a semisimplicial set and a simplicial complex by the same letter, then they are related by $X_\circ^{\ord} = X_\bullet$. By $X_p$ we will always mean the set of $p$--simplices of the semisimplicial set $X_\bullet$, not the set of $p$--simplices of the simplicial complex $X_\circ$, which is then in bijection with $X_p/ \mathfrak S_{p+1}$. 
\end{convention}

 %%%%%%%%%%%%% START PURPLE

The following theorem is well known; see for example Kupers--Miller \cite[Lemma 3.16]{KM2} and Randal-Williams--Wahl \cite[Proposition 2.14]{RWW}. 

\begin{theorem} \label{complextosemi} \qquad 
\begin{enumerate} 
\item[(a)] If $X_\bullet = X^{\ord}_\circ$ is $n$--connected, then so is $X_\circ$. 
\item[(b)] If $X_\circ$ is weakly Cohen--Macaulay of dimension $n$, then $X_\bullet = X^{\ord}_\circ$ is $(n-1)$--connected. 
\end{enumerate}
\end{theorem}
 %%%%%%%%%%%%%%%%%%%%%% END PURPLE

These results will allow us to pass between high connectivity results for semisimplicial sets and simplicial complexes. The following definition and theorem is a useful tool for proving simplicial complexes are highly connected.

\begin{definition} \label{defjoin}
A map of simplicial complexes $\pi\colon Y_\circ \to X_\circ$ is said to exhibit  $Y_\circ $ as a \emph{join complex} over $X_\circ$ if it satisfies all of the following:  
\begin{itemize}
  \item $\pi$ is surjective
  \item $\pi$ is simplex-wise injective, that is, the  induced map on the geometric realizations is injective when restricted to a (closed) simplex
  \item a collection of vertices $(y_0,\ldots, y_p)$ spans a simplex of $Y$ whenever there exists simplices $\theta_0, \ldots, \theta_p$ such that for all $i$, $y_i$ is a vertex of  $\theta_i$ and the simplex $\pi(\theta_i)$ has vertices $\pi(y_0), \ldots, \pi(y_p)$.
\end{itemize}
\end{definition}

%This definition is illustrated in  \autoref{DefnJoinComplex}. 
%
%
%\begin{figure}[!ht]    \centering
%\labellist
% \hair 0pt
%\pinlabel {\color{blue} \scriptsize $\theta$} [c] at  37 34
%\pinlabel {\color{blue} \scriptsize $\theta_0$} [c] at  20 47
%\pinlabel {\color{blue} \scriptsize $\theta_1$} [c] at  20 20 
%\pinlabel {\color{black} \scriptsize  $\pi$} [c] at 56 39
%\pinlabel {\scriptsize  $Y_{\circ}$} [c] at 17 0
%\pinlabel {\scriptsize $X_{\circ}$} [c] at  77 2
%\pinlabel {\scriptsize $y_0$} [c] at  31 57
%\pinlabel {\scriptsize $y_1$} [c] at  31 8
%\pinlabel {\scriptsize $y_2$} [c] at  1 57
%\pinlabel {\scriptsize $y_3$} [c] at  1 8
%\pinlabel {\scriptsize $\pi(y_0)=\pi(y_2)$} [l] at  81 52
%\pinlabel {\scriptsize $\pi(y_1)=\pi(y_3)$} [l] at  81 14
%\endlabellist
%\begin{center}\scalebox{1.3}{\includegraphics{DefnJoinComplex}}\end{center}
%\caption{The map $\pi$ does not exhibit $Y_{\circ}$ as a join complex over $X_{\circ}$ unless $\theta$ is a simplex of $Y_{\circ}$.}
%\label{DefnJoinComplex}
%\end{figure}  
%
%
%\begin{definition}
%Let $X_\circ$ a simplicial complex and for each simplex $\sigma$ and vertex $x \in \sigma$, let $L_x(\sigma)$ be a set. The collection of sets $\{L_x(\sigma)\}$ is called a labeling system if $L_x(\sigma) \subseteq L_x(\tau)$ whenever $x \in \tau \subseteq \sigma$.  
%\end{definition}

The following result is due to Hatcher--Wahl \cite{hatcherwahl}.

\begin{theorem}[Hatcher--Wahl {\cite[Theorem 3.6]{hatcherwahl}}] \label{labelsystemconn}
Let $\pi\colon Y_\circ \to X_\circ$ be a map of simplicial complexes exhibiting $Y_\circ$ as a join complex over $X_\circ$. Assume $X_\circ$ is wCM of dimension $n$.  Further assume that for all $p$-simplices $\tau$ of $Y_\circ$, the image of the link $\pi(\Lk^Y(\tau))$ is wCM of dimension $(n-p-2)$. Then $Y_{\circ}$ is $\displaystyle \left \lfloor \frac{n-2}{2}\right \rfloor$--connected. 
\end{theorem}

\begin{definition}
Let $X_\circ$ be a simplicial complex and  $X_\bullet = X^{\ord}_\circ$. Let $\widetilde\sigma \in X_p$ and $\sigma\subseteq X_\circ$ the corresponding simplex. Then define the \emph{link} 
\[ \Lk_\bullet^X \widetilde \sigma = (\Lk_\circ^X \sigma)^{\ord}\]
as a sub-semisimplicial set of $X_\bullet$.
\end{definition}

\subsection{Algebraic preliminaries} \label{SecAlgebraReview}

In this subsection, we recall some basic facts and definitions concerning free modules over PIDs as well as symplectic structures.  

Let $R$ be a PID. For $S$ a set, throughout this paper we write $R[S]$ to denote the free $R$--module with basis $S$. Notably, this does \emph{not} denote the polynomial algebra with variables $S$. 

Given a submodule $W$ of a free $R$--module $V$, we say $W$ is \emph{splittable} or \emph{has a complement} if there exists a submodule $U$ with $V=W \oplus U$. Given a submodule $W$ of a free module $V$, let $\sat(W)$ denote the intersection of all splittable submodules of $V$ which contain $W$, equivalently,   $\sat(W)$ is the preimage of the torsion submodule of $V/W$ under the quotient map $V \to V/W$. We call $\sat(W)$ the \emph{saturation} of $W$. We write $\rk(W)$ for the rank of a free module. The following proposition collects some elementary facts concerning free modules over PIDs and their submodules; see for example Kaplansky \cite{Kaplansky}.

\begin{proposition} \label{modulefacts}
Let $R$ be a PID and $A,B,C$ submodules of a finitely generated free $R$--module $V$.
 
\begin{enumerate}
  \item  If $A$ and $B$ have complements in $V$, then so does $A \cap B$.
  \item Let $B$ have a complement. Then $A$ has a complement containing $B$ if and only if there is a submodule $D$ with $V = A \oplus B \oplus D$. \label{containsimpliessum}
  \item $\rk(A)=\rk(\sat(A))$.  
  \item $\Big( (A \cap C)\oplus (B \cap C) \Big)\subseteq (A + B)\cap C$ but equality does not hold in general. 
  \item \label{sub} If $V=A \oplus B$ and $C \supseteq A$, then $C=A \oplus (B \cap C)$.
  \label{iv}
  \item If $V=A \oplus B$ and $C \subseteq A$, then $C=A \cap (B \oplus C)$.
  \label{v}
  \item $A$ has a complement in $V$ if and only if $V/A$ is torsion free.
%  { \color{Red} %%%%%%%%%%%%%%%%%%%%%%%%%%%%%%%%%%%%% BEGIN RED
%  \item If $A \subseteq B$ and $A$ has a complement in $V$, then $A$ has complement in $B$.  
%  }          %%%%%%%%%%%%%%%%%%%%%%%%%%%%%%%%%%%%% END RED
\end{enumerate} 
\end{proposition}

   %%%%%%%%%%%%%%%%%%%%%%%%%%%%%%%%%%%%% BEGIN PURPLE
%Notably,  \autoref{modulefacts} \autoref{sub} implies that if $A \subseteq B$ and $A$ has a complement in $V$, then $A$ has complement in $B$.  
            %%%%%%%%%%%%%%%%%%%%%%%%%%%%%%%%%%%%% END PURPLE

 Again let $R$ be a PID and $V$ a finitely generated free $R$--module. Recall that a \emph{symplectic form} on $V$ is a perfect bilinear pairing $\langle\,,\,\rangle \colon V \times V \to R$ such that $\langle v,v\rangle =0$ for all $v\in V$. %A submodule $W \subseteq V$ is called \emph{isotropic} if $\langle w_1,w_2\rangle = 0$ for all $w_1,w_2 \in W$. 
 For a submodule $W \subseteq V$, we write $W^\perp$ to denote its symplectic complement
\[ W^\perp = \{ v \in V \mid \langle v, w\rangle = 0 \text{ for all $w \in W$}\}.\]
 % {\bf \color{Purple} %%%%%%%%%%%%%%%%%%%%%%%%%%%%%%%%%%%%% BEGIN PURPLE
Recall that a submodule $W$ of a symplectic $R$--module $V$ is called \emph{symplectic} if  the form restricts to a symplectic form on $W$.
%and \emph{isotropic}  if $W \subseteq W^{\perp}$. 
We note that, when $R$ is not a field, the condition that $W$ is symplectic may be stronger than the condition that $W \cap W^{\perp} = 0$. 

 % }          %%%%%%%%%%%%%%%%%%%%%%%%%%%%%%%%%%%%% END PURPLE
The following proposition is an elementary fact about free symplectic modules over PIDs and their submodules (see for example Knus \cite[4.1.2]{Knus}).

\begin{proposition} \label{SymplecticBasis} 
Let $R$ be a PID and $V$ a finitely generated free symplectic module over $R$. Then the rank of $V$ is even and there is a basis $v_1,w_1, \dots, v_n,w_n$ of $V$ such that
  \[ \langle v_i,v_j \rangle = \langle w_i, w_j\rangle = 0, \quad \langle v_i, w_j\rangle = \delta_{ij} \qquad \text{for all $i, j=1, \ldots, n$}.\]
  \end{proposition}
 
 %\begin{proposition} \label{SymplecticFacts}
%Let $R$ be a PID and let $A,B,C$ be submodules of a finitely generated free symplectic module $V$.
%\begin{enumerate}
%  \item \label{SymplecticBasis}  The rank of $V$ is even and there is a basis $v_1,w_1, \dots, v_n,w_n$ of $V$ such that
%  \[ \langle v_i,v_j \rangle = \langle w_i, w_j\rangle = 0, \quad \langle v_i, w_j\rangle = \delta_{ij} \qquad \text{for all $i, j=1, \ldots, n$}.\]
%  \item \label{ii} If $A$ is a direct summand in $V$, every maximal symplectic submodule $U$ of $A$ has the same rank and is splittable. In particular
%  \[ A = U \oplus (U^\perp  \cap A),\]
%  and $U^\perp \cap A$ is isotropic.
%\end{enumerate} 
% \end{proposition}
%%%%%%%%%%%%%%%%%%%%%%%%%%%%Part (ii) was not true as written. 

%\begin{proof}
 % {\bf \color{Purple} %%%%%%%%%%%%%%%%%%%%%%%%%%%%%%%%%%%%% BEGIN PURPLE
 % Let $\iota$ denote the map 
%  \begin{align*} V & \longrightarrow V^* \\ 
%   v & \longmapsto [ w \mapsto \langle v, w \rangle] 
%  \end{align*} 
%  Recall that the assumption that the form $\langle -,- \rangle$ is \emph{perfect} is the statement that $\iota$ is an isomorphism. 
%The proposition  is a standard result; see for example McDuff--Salamon \cite[proof of Theorem 2.1.3]{McDuffSalamon}. McDuff--Salamon state the result for the case $R=\R$, but (when we replace the hypothesis that the form is non-degenerate with the assumption that it is perfect) their proof holds for $R$ a PID. 
%  }          %%%%%%%%%%%%%%%%%%%%%%%%%%%%%%%%%%%%% END PURPLE
%\end{proof}

\subsection{Generalized partial basis complexes} \label{SecPB} 

In this subsection, we recall the definition of several simplicial complexes and semisimplicial sets involving partial bases and complements.  These connectivity results will imply that $\VIC(R)$--modules satisfying a polynomial condition exhibit central stability. Throughout the subsection we let $R$ be a PID and $V$ be a finite-rank free $R$--module. 

\begin{definition}
A \emph{partial basis} of a free module $V$ is a  linearly independent set $\{v_0,\ldots, v_p\} \subseteq V$ such that there is a free (possibly zero) submodule $C$ with $\spn(v_0,\ldots,v_p) \oplus C = V$. Such a set  $\{v_0,\ldots, v_p\}$ is also called \emph{unimodular}, and the submodule $C$ is called a \emph{complement} for the partial basis. An \emph{ordered partial basis} is a partial basis with a choice of bijection with a set of the form $\{0,\ldots,p\}$. 
\end{definition}

\begin{definition} 

\label{defPB}
For $p \ge 0$, let $\PB_p(V)$ be the set of ordered partial bases of size $p+1$. For $0 \le i \le p$, there are maps $d_i\colon \PB_p(V) \to \PB_{p-1}(V)$ given by forgetting the $i$th basis element. With these maps, the sets  $\PB_p(V)$ assemble into a semisimplicial set $\PB_\bullet(V)$. Let $\PB_\circ(V)$ be the unique simplicial complex with $\PB_\circ^{\ord}(V)=\PB_\bullet(V)$

 \end{definition} 

The existence of $\PB_\circ(V)$ and all other simplicial complexes considered in this subsection follow easily from the work of Randal-Williams--Wahl \cite[Definition 2.8 and proof of Lemma 5.10]{RWW}. 

 The link of a simplex $\sigma=(v_0, v_1, \ldots,v_p) \in \PB_p(V)$ is the subcomplex of ordered partial bases $(u_0, u_1, \ldots, u_q)$ such that $\{v_0, v_1, \ldots,v_p, u_0, u_1, \ldots, u_q\}$ is a partial basis of $V$. Notably, this link only depends on the submodule $W =\spn(v_0,\ldots,v_p)$. By abuse of notation, we will often denote $\Lk_\bullet(\sigma)$ by $\Lk_\bullet(W)$. By convention if $W=0$, we let $\Lk_\bullet(W)$ be the entire complex of ordered partial bases $\PB_\bullet(V)$. If $U\subseteq V$ is a splittable submodule, then there is a canonical inclusion of $\PB_\bullet(U) \subseteq \PB_\bullet(V)$. 

\begin{definition} 
\label{defI} 
Let $U,W\subseteq V$ be splittable submodules. Define
\[ \PB_\bullet(U,W) = \PB_\bullet(U) \cap \Lk_\bullet(W)\subseteq \PB_\bullet(V)\]
and let $\PB_\circ(U,W)$ be the unique simplicial complex with $\PB_\circ^{\ord}(U,W)=\PB_\bullet(U,W)$

 \end{definition} 

%It is straightforward to check that $\PB_\circ(U,W)$ exists.
Concretely, $\PB_\bullet(U,W)$ is the sub-semisimplicial set of $\PB_\bullet(U)$ consisting of ordered nonempty partial bases of $U$ contained in a complement of $W$. The complex $\PB_\bullet(U,W)$ depends only on the submodule sat$(U+W)$ and not on $V$. We note that \[\PB_\bullet(U,W) = \Lk^{\PB(U)}_\bullet(W) \qquad \text{if $W \subseteq U$},\]  and in particular that 
\[  \PB_{\bullet}(V,V)=\varnothing \qquad \text{and} \qquad  \PB_{\bullet}(V,0)=\PB_{\bullet}(V).\] 
More generally, whenever $W$ is contained in any complement of $U$, then $\PB(U,W) \cong \PB(U)$.

\begin{remark}
For vector spaces,  \[\PB_\bullet(U,W) = \PB_\bullet(U,W\cap U) = \Lk^{\PB(U)}_\bullet(W \cap U) \qquad \text{($R$ a field).}\]
We caution, however, that this identification does not hold in general. For example, when $R=\Z$ and $V=\Z^3$, consider the submodules 
\[U=\spn \left(\begin{bmatrix}1 \\ 1 \\ 0 \end{bmatrix} \right) \qquad \text{and} \qquad W=\spn \left(\begin{bmatrix}1 \\ -1 \\ 0 \end{bmatrix} \right).\]
Since the determinant of the matrix \[ \begin{bmatrix}1 & 1 & a \\ 1 & -1 & b \\ 0 & 0 & c \end{bmatrix}\] is a multiple of $2$ for any $a,b,c \in \Z$, there is no basis for $V$ that contains both a basis for $U$ and a basis for $W$, and $\PB_\bullet(U,W)$ is empty.  In contrast, the complex $\PB_\bullet(U,W\cap U) =  \PB_\bullet(U,0) =  \PB_\bullet(U)$ is nonempty. 
\end{remark}

%\todo{ADD LEMMA -- can still realize $\PB_\bullet(U,W)$ as some link in $\PB_\bullet(U)$ -- necessary for using Maazen results.  } %Peter: Don't think this lemma is true.

%\begin{definition} 
%\label{defI} 
%Let $U,Q\subseteq V$ be splittable subspaces. Pick an ordered basis $(q_0,\ldots,q_j)$ of $Q$. Define $\PB_p(U,Q) \subseteq \PB_p(U)$ is the set of ordered partial bases $(u_0,\ldots,u_p)$ of $U$ such that $(u_0,\ldots,u_p,q_0,\ldots,q_j)$ is an ordered partial basis of $V$.
%Let $\PB_\circ(U,Q)$ to be the simplicial complex formed by quotienting by the symmetric group actions.
%\end{definition}

%Note that the definition of $\PB_\circ(U,Q)$ does not depend on the choice of basis of $Q$. 

%If $R$ is a field and $U \cap Q=0$, then $\PB_\circ(U,Q)=\PB_\circ(U)$. However, for general rings this is not always the case. For example, if $R=\Z$, $V=\Z^2$, $U=\spn \left(\begin{bmatrix}
%1 \\ 0 \end{bmatrix} \right)$ and $Q=\spn \left(\begin{bmatrix}
%1 \\ 2 \end{bmatrix} \right)$, then $\PB_\circ(U,Q)$ is empty while $\PB_\circ(U)$ is not empty. We now define some more semisimplicial sets and simplicial complexes. These will involve the data of complements of partial bases. 

We next define a variation of the semisimplicial set $\PB_{\bullet}(V)$ consisting of ordered partial bases with distinguished choices of complements. 

\begin{definition} 
\label{defPBC}
Let $\PBC_p(V)$ be the set of ordered partial bases $(v_0,\dots,v_p)$ of $V$ as well as a choice of complement $C$ such that
\[ C\oplus \spn(v_0,\dots, v_p)  = V.\] Let $d_i\colon \PBC_p(V) \to \PBC_{p-1}(V)$ be given by the formula \[d_i(v_0,\ldots,v_p,C)=(v_0,\ldots,\hat v_i, \ldots, v_p,C \oplus R v_i).\] Here the hat indicates omission. These sets assemble to form a semisimplicial set $\PBC_\bullet(V)$. Let $\PBC_\circ(V)$ be the unique simplicial complex with $\PBC_\circ^{\ord}(V)=\PBC_\bullet(V)$

 \end{definition} 

%It is straightforward to check that $\PBC_\circ(V)$ exists given the following proposition. 

\begin{proposition} \label{alternatedescription}
Let $\partial_i:\PBC_p(V) \m  \PBC_0(V)$ be the map induced by $\{i\} \hookrightarrow \{0,\ldots, p\}$ and let $\partial:\PBC_p(V) \m \PBC_0(V)^{p+1}$ be the product of these maps. The map $\partial$ is injective. If $R$ is a PID, then \[\big( (v_0,C_0), \ldots, (v_p,C_p) \big)\] is in the image of $\partial$ if and only if $v_i \in C_j$ for all $i \neq j$. 
\end{proposition}

\begin{proof}
Consider  $(v_0, \dots, v_p,C) \in \PBC_p(V)$. We have that \[\partial_i(v_0, \dots, v_p,C) =(v_j, C\oplus \spn(v_0\dots,\hat v_j, \dots, v_p)) \in \PBC_0(V).\] The ``only if'' portion of the claim follows from this description. 

Conversely, let $\big((v_0,C_0), \ldots, (v_p,C_p) )\in (\PBC_{0}(V))^{p+1}$ such that $v_i \in C_j$ for $i\neq j$. We will prove that
\[ (v_0, \dots, v_p, C) \in \PBC_p(V), \qquad \text{where} \qquad C = \bigcap_{i=0}^p C_i \ . \] This will imply the ``if'' portion of the claim as then  \[\partial(v_0, \dots, v_p, C) = \big((v_0,C_0), \ldots, (v_p,C_p) \big). \]  Define
\[ D_j = \bigcap_{i=0}^j C_i,\]
and we will prove by induction that
\[ \spn(v_0, \dots, v_j) \oplus D_j = V.\]
The base case is the statement that $\spn(v_0) \oplus C_0=V$. 
%By \autoref{modulefacts} \autoref{v} and the inductive hypothesis, 
%\begin{align*} &\spn(v_0, \dots, v_j) \cap D_j  \\
%& = \Big( \spn(v_0, \dots, v_{j-1})\oplus \spn(v_j)\Big) \cap D_{j-1} \cap C_j \\
%& =  \spn(v_j) \cap C_j \\
%&= \{ 0\}.
%\end{align*}
Since $\spn(v_0, \ldots, v_{j-1}) \subseteq C_j$ by assumption, \autoref{modulefacts} \autoref{iv} and the inductive hypothesis imply that
\[ \spn(v_0, \dots, v_{j-1}) \oplus (D_{j-1} \cap C_j) = C_j.\]
By taking the direct sum of both sides of this equation with $\spn(v_j)$, we conclude the inductive step. Applying this result when $j=p$ yields the desired decomposition
\[ \spn(v_0, \dots, v_p) \oplus C = \spn(v_0, \dots, v_p) \oplus D_p = V. \] This establishes the ``if'' portion of the claim. 

Injectivity of $\partial$ follows from the fact that if \[\partial(v_0, \dots, v_p, C) = \big((v_0,C_0), \ldots, (v_p,C_p) \big), \] then \[C = \bigcap_{i=0}^p C_i. \qedhere \]

\end{proof}
%We remark in particular that the $j$th vertex of a simplex $(v_0, \dots, v_p,C)$ in $\PBC_\bullet(V)$ is
%\[ (v_j,C_j) = (v_j, \spn(v_0, \dots, \hat v_j, \dots, v_p) \oplus C).\]

It will be convenient for us to realize the complexes $\PBC_\bullet(V)$ and their links as special cases of the following more general construction. 
In the following, we will often identify an ordered partial basis $(v_0,\dots,v_p)$ of $V$ with a split $R$--linear monomorphism $f \colon R^{p+1}\to V$.
\begin{definition}
Let $U,W\subseteq V$ be splittable submodules. Define the sub-semisimplicial set $\PBC_\bullet(V,U,W) \subseteq \PBC_\bullet(V)$ by
\[ \PBC_p(V,U,W) = \{ (f,C)\in \PBC_p(V) \mid \im f \subseteq U, \; W \subseteq C\}.\]
Let $\PBC_\circ(V,U,W)$ be the unique simplicial complex with $\PBC_\circ^{\ord}(V,U,W)=\PB_\bullet(V,U,W)$
 \end{definition} 

%Using \autoref{alternatedescription}, it is straightforward to check that $\PBC_\circ(V,U,W)$ exists. 
Note that \[ \PBC_\bullet(V) \cong   \PBC_\bullet(V,V,0).\] 
%There is a natural bijection $\PBC_p(V)\cong\Hom^{\VIC}(R^{p+1},V)$. In other words, the data of a complemented partial basis is the same data as an injective homomorphisms $g\colon R^{p+1} \to V$ and free submodules $D \subseteq V$ with $\im(g) \oplus D=V$.
Given a simplex $\sigma = (f,C) \in \PBC_p(V)$, we can identify its link
\[ \Lk_\bullet(\sigma) = \PBC_\bullet(V,C,\im f).\]
This link is isomorphic to $\PBC_\bullet(C)$, which we will see is a special case of \autoref{simplify} below. 
%The embedding of $\PBC_\bullet(C)$ into $\PBC_\bullet(V)$ depends on the image of $f$:
%\begin{align*}
%\PBC_p(C) &\longrightarrow \Lk_\bullet(\sigma)\subseteq \PBC_p(V)\\
%(c_0,\dots, c_p,D) &\longmapsto (c_0,\dots, c_p,\im f\oplus D)
%\end{align*}
In general $\PBC_\bullet(V,U,W)$ is \emph{not} isomorphic to $\PBC_\bullet(U)$, but the following lemma gives a more general picture. 

\begin{lemma}\label{simplify}
Let $R$ be a PID, $U,W\subseteq V$ be splittable, and $A\oplus B = V$ such that $U \subseteq A$ and $B \subseteq W$. Then
\begin{align*} \PBC_\bullet(V,U,W) &\overset{\cong}\longrightarrow \PBC_\bullet(A,U,W\cap A)\\
(f,C) & \longmapsto (f, C\cap A)
 \end{align*}
is an isomorphism. In particular, if $U \oplus W = V$, then $ \PBC_\bullet(V,U,W) \cong \PBC_\bullet(U)$.  
\end{lemma}

\begin{proof}
Let $(f,C)\in \PBC_\bullet(V,U,W)$. \autoref{modulefacts}  \autoref{iv} implies \[ \im f \oplus (C\cap A) = A\] and so we conclude that $(f,C\cap A)\in \PBC_\bullet(A,U,W\cap A)$.  If $(g,D)\in \PBC_\bullet(A,U,W\cap A)$, then $(g, D\oplus B) \in \PBC_\bullet(V,U,W)$. These two maps are inverses because
\[ (C\cap A) \oplus B = C\]
by  \autoref{modulefacts}  \autoref{iv}, and by  \autoref{modulefacts}  \autoref{v}, 
\[ (D \oplus B) \cap A = D. \qedhere\]
%When  $U \oplus W = V$, taking $A=U$ and $B=W$ gives the special case \[ \PBC_\bullet(V,U,W)  \cong   \PBC_\bullet(U,U,0) =   \PBC_\bullet(U). \qedhere \] 
\end{proof}

%Let $e_0,\ldots, e_p$ be the standard basis of $R^{p+1}$. Let $g\colon R^{p+1} \to V$ be an element of $\PB_p(V)$. This is same data as an ordered partial basis $(g(e_0),\ldots,g(e_p))$ of $V$. Similarly, a $p+1$ simplex of  $\PB_\circ(V)$ is an unorderd partial basis of $V$. 

%Note that ordered partial basis of $V$ of size $d$ is the same data as an $R$--linear injection $f:R^d \m V$ with splittable image. We will often use this perspective on partial bases. From now on, let $\alpha=(V,X,Y,Z,W)$ with $X,Y,Z$ and $W$ free submodules of $V$ with $X \oplus Y=V$ and $Z \oplus W=V$.
%
%
%\begin{definition}  \label{PBCa}
%Let $\PBC^{\alpha}_p$ be the set of pairs $(f,C)$ with $f\colon R^{p+1} \to X \cap Z$ an $R$--linear injection, $C \subseteq V$ a submodule containing $Y$ and $W$ such that $V=\im(f) \oplus C$. Viewing $f$ as an ordered partial basis of $X\cap Z$, the sets $\PBC_\bullet^\alpha$ assemble to form a semisimplicial set with face maps given by removing vectors from the partial basis and adding them to the subspace $C$ as in  \autoref{defPBC}. Define $\PBC^\alpha_\circ$ to be the simplicial complex formed by quotienting by the symmetric group actions. 
%%old
%\end{definition}

%For $W=Y=0$ and $X=Z=V$, $\PBC^{\alpha}_\bullet = \PBC_\bullet(V)$.  Let $t$ be $\rk(X) + \rk(Z)-\rk(V)$ and let $n$ be the rank of $V$. Generically, $t$ is the rank of $X \cap Z$. Let $R$ be a ring with  stable rank $s$. The main theorem of this subsection is the following. 

The main theorem of this subsection is  \autoref{PBCsemiconn}. To state this theorem, we will use Bass' stable range condition for rings. 

%\begin{definition} \label{DefnStableRank} Given a positive integer $s$, a ring $R$ has \emph{stable rank $s$} if whenever \[a_0R+ a_1R + \cdots  + a_mR = R, \qquad a_i \in R,\] with $m \ge s$, there exist elements $x_1, x_2, \ldots, x_m \in R$ such that \[(a_1 +a_0x_1)R + \ldots + (a_m + a_0x_m)R = R.\]  \end{definition}

\begin{definition}[Bass {\cite[Section 4 ``Definition"]{BassKTheory} and \cite[Definition 3.1]{Bass}}] \label{DefnStableRank} Let $s$ be a positive integer. A ring $R$ has \emph{stable rank $s$} if $s$ is the smallest positive integer $m$ for which the following \emph{Condition ($B_m$)} holds: 
 whenever 
\[a_0R+ a_1R + \cdots  + a_mR = R, \qquad a_i \in R,\]
there exist elements $x_1, x_2, \ldots, x_m \in R$ such that 
\[(a_1 +a_0x_1)R + \ldots + (a_m + a_0x_m)R = R.\]
  \end{definition}

%The condition was first formulated by Bass to characterize when matrices in GL$_{m+1}(R)$ can be row-reduced to matrices in the image of GL$_{m}(R)$.
%, a finite direct sum of fields, or the ring $\Z/m\Z$ for $m \in \Z$,

 Note that the indexing in \autoref{DefnStableRank} differs from Bass' convention in his work \cite{Bass};  he called it the stable range condition $SR_{s+1}$.  Bass (\cite[Section 4 ``Examples"]{BassKTheory} and \cite[Proposition 3.4 (a) and Theorem 3.5]{Bass})  states that semi-local rings have stable rank $1$, and Dedekind domains have stable rank at most $2$. So, for example, $R$ has stable rank $1$ if $R$ is a field and $R$ has stable rank at most $2$ when $R$ is a PID. Since the next result concerns PIDs, the stable rank $s$ must be equal to 1 or 2. 
 
 The remainder of the subsection serves to prove the following theorem. 
 %, or when $R$ is the ring of integers of a number field. The integers $R=\Z$ have stable rank $2$. More generally, a commutative Noetherian ring of Krull dimension $d$ has stable rank at most $(d+1)$.  %\cite[Corollary 6.5]{Bass} is another reference for the semi-local ring statement. 

\begin{theorem} \label{PBCsemiconn} 
Let $R$ be a PID. Let $s$ be the stable rank of $R$. If $U$ and $W$ are splittable submodules of the finite-rank free  $R$--module $V$, then $ \PBC_\bullet(V,U,W)$
is $\displaystyle \left \lfloor \frac{\rk U - \rk W-s-2}{2}\right \rfloor$--connected.
\end{theorem}

 \autoref{PBCsemiconn} is a partial generalization of the result of  \cite[Lemma 5.10]{RWW} of Randal-Williams--Wahl that if $R$ has stable rank $s$ then $\PBC_\bullet(V)$ is $\displaystyle \left \lfloor \frac{\rk V-s-2}{2}\right \rfloor$--connected. 
%\begin{theorem}[Randal-Williams--Wahl {\cite[Lemma 5.10]{RWW}}]  \label{PBCVconn} If $R$ is a ring with stable rank $s$, then $\PBC_\bullet(V)$ is $\displaystyle \left(\frac{\rk V-s-2}{2}\right)$--connected. \end{theorem}
Following their proof, we will prove  \autoref{PBCsemiconn} by comparing $\PBC_\circ(V,U,W)$ with its image in $\PB_\circ(V)$. This will let us use a result of van der Kallen {\cite[Theorem 2.6 (i)--(ii)]{V}}: Assume $R$ is a ring with stable rank $s$. Let $U$ and $W$ be splittable subspaces of $V$. Then 
$\PB_\bullet(U,W)$
 is $(\rk U-\rk W-1-s)$--connected. In particular, $\PB_\bullet(U)$ is $(\rk U-1-s)$--connected.

Note that van der Kallen's result is stated in terms of the poset whose underlying set is \[ \bigsqcup_p \PB_p(U,W) \] and whose order is induced by the face maps. The geometric realization of this poset is isomorphic to the barycentric subdivision of the geometric realization of $\PB_\bullet(U,W)$ (see for example Randal-Williams--Wahl \cite[Proof of Lemma 5.10]{RWW}) and hence the connectivity results apply. We will use this same identification in the proof of \autoref{USpRH4}.

%This gives the following corollary. 

\begin{proposition}\label{Vancor} Let $R$ be a ring with stable rank $s$, and let $U$ and $W$ be splittable subspaces of the finite-rank free $R$--module $V$.
The simplicial complex  $\PB_\circ(U,W)$
is wCM of dimension $(\rk U-\rk W-s)$. In particular, $\PB_\circ(U)$ is wCM of dimension $(\rk U-s)$.  
\end{proposition}

\begin{proof}
By  \cite[Theorem 2.6 (ii)]{V} and  \autoref{complextosemi} (a), the complex $\PB_\circ(U,W)$ is $(\rk U-\rk W-1-s)$--connected. The link of $\{v_0,\ldots, v_p\}$ in $ \PB_\circ(U,W) = \PB_\circ(U)\cap \Lk_\circ(W)$ is isomorphic to 
\[\PB_\circ(U,W\oplus\spn(v_0,\ldots, v_p)  ) = \PB_\circ(U) \cap\Lk_\circ(W\oplus\spn(v_0,\ldots, v_p)  )\]
 and so the links are $(\rk U-\rk W-s-2-p)$--connected as required.
\end{proof}

To prove $\PBC_\bullet(V,U,W)$ is highly connected, we will show that $\PBC_\circ(V,U,W)$ is wCM. To do this, we need the following lemma.

\begin{lemma} \label{linksPBC}
Let $R$ be a PID and $U,W,X,Y,Z$ be splittable submodules of a finite-rank free $R$--module $V$.
\begin{enumerate}
\item \label{item:intersectionPBC}  $\PBC_\bullet(V,X,Y)\cap \PBC_\bullet(V,Z,W) = \PBC_\bullet(V,X\cap Z, \sat(Y+W))$.
\item \label{item:linkPBC} Any simplex $\sigma = (f,C) \in \PBC_\bullet(V,U,W)$ has link $\Lk_\bullet(\sigma) = \PBC_\bullet(C,U\cap C, W)$.
\end{enumerate}
\end{lemma}

\begin{proof}
\begin{enumerate}
\item Both sides of the equation describe the following semisimplicial set:
\begin{align*} & \{ (f,C)\in \PBC_\bullet(V) \mid \im f \subseteq X, \im f\subseteq Z, Y\subseteq C, W\subseteq C\} \\ & = \{ (f,C)\in \PBC_\bullet(V) \mid \im f \subseteq X\cap  Z, \sat(Y+W)\subseteq C\} \end{align*}
\item 
Every simplex in $\PBC_\bullet(V,U,W)$ contains every simplex in  $\PBC_\bullet(V)$ spanned by vertices in $\PBC_\bullet(V,U,W)$ (in other words, the inclusion $\PBC_\bullet(V,U,W) \hookrightarrow \PBC_\bullet(V)$ is \emph{full}). Hence, for $\sigma \in \PBC_\bullet(V)$, 
%Because $\PBC_\bullet(V,U,W)\subseteq \PBC_\bullet(V)$ is full,
\begin{align*} \Lk_\bullet^{\PBC_\bullet(V,U,W)}(\sigma) &= \Lk_\bullet^{\PBC_\bullet(V)}(\sigma) \cap \PBC_\bullet(V,U,W) \\ 
&= \PBC_\bullet(V,C,\im f) \cap \PBC_\bullet(V,U,W) \\ 
&= \PBC_\bullet(V, U\cap C, W \oplus \im f). 
\end{align*}
The last step uses \autoref{item:intersectionPBC} and the observation that, since $W$ is contained in the complement $C$ of $\im f$,  \[\sat(W+\im f) = W \oplus \im f.\] By applying  \autoref{simplify} with $A=C$ and $B=\im f$, we find
\[  \PBC_\bullet(V, U\cap C, W \oplus \im f) =  \PBC_\bullet\Big(C, U\cap C, (W \oplus \im f)\cap C\Big).\] 
 Then \autoref{modulefacts}  \autoref{v} implies that 
\[ (W \oplus \im f) \cap C = W\]
and the result follows. \qedhere
\end{enumerate} 
%Let $\tau$ be in the link of $\sigma$ and let $(g,D)\in \PBC_\bullet^\alpha$ be a preimage of $\tau$. Then the join $\sigma * \tau $ is a simplex in $\PBC_\circ^{\alpha}$. Because of how the face maps work, necessarily $\im(f) \subseteq X \cap Z \cap D$ and $\im(g) \subseteq X \cap Z \cap C$ are direct summands and a preimage of $\sigma*\tau$ in $\PBC_\bullet^\alpha$ has the form $(h,C \cap D )$ with $\im(h) = \im(f) \oplus \im(g)$. 
% 
%We define a map from $\Lk(\sigma)$ to $\PBC_\circ^\beta$ by sending $(g,D)$ to $(g,C \cap D)$. For well-definedness, we note that we have already observed that $\im(g) \subseteq X \cap Z \cap C$. Since $C$ and $D$ each contain $Y$ and $W$, so does $C \cap D$. We must also check that $C = \im(g) \oplus (C \cap D) $. This follows from  \autoref{iv} of  \autoref{modulefacts}.  
% 
%We claim that $(g,D) \mapsto (g, D\oplus\im f)$ induces an inverse map $\PBC_\circ^\beta \to \Lk(\sigma)$. This map is easily checked to be well defined. To see that it is an inverse, use  \autoref{iv} and \autoref{v} of   \autoref{modulefacts}. 
\end{proof}

%\begin{theorem} \label{PBconn}
%The simplicial complex $\PB(V)_\circ$ is wCM of dimension $n-s$. 
%\end{theorem}

%We now need an alternate description of simplices in $\PBC_\circ(U,W)$. Its proof is elementary and uses   \autoref{containsimpliessum} of  \autoref{modulefacts}.

\begin{proposition} \label{CorPBCwCM}
If $R$ is a PID of stable rank $s$, the complex $\PBC_\circ(V,U,W)$ is wCM of dimension $\displaystyle \left \lfloor \frac{\rk U-\rk W-s}{2}\right \rfloor$. 
\end{proposition}

\begin{proof}

%To show is wCM of dimension $\displaystyle \left \lfloor \frac{\rk U-\rk W-s}{2}\right \rfloor$,  it suffices to check that the simplicial complex $\PBC_\circ(V,U,W)$ is $\displaystyle \left \lfloor \frac{\rk U -\rk W-s-2}{2}\right \rfloor$--connected.

%The proposition will follow by applying  \autoref{labelsystemconn} to the map $\theta\colon \PBC_\circ(U,W) \to \PB_\circ(U,W)$.

We begin by checking that the map 
\begin{align*}
\theta\colon \PBC_\circ(V,U,W) &\longrightarrow \PB_\circ(U,W) \\
 \{(v_0,C_0), \dots, (v_p,C_p)\} & \longmapsto \{v_0, \dots, v_p\}
\end{align*} exhibits $ \PBC_\circ(V,U,W)$ as a join complex over $\PB_\circ(U,W)$. Simplex-wise injectivity is clear. We next check  surjectivity. Let $\{v_0,\dots,v_p\}$ be a $p$--simplex in $\PB_\circ(U,W)$. Because $\{v_0,\dots,v_p\} \in \Lk_\circ(W)$, there is a submodule of $D \subseteq V$ such that
\[ D \oplus \spn(v_0,\dots,v_p)\oplus W = V.\] 
Then the simplex $\{(v_0,C_0), \dots, (v_p,C_p)\}$ with $C_i = \spn(v_0, \dots,\hat v_i, \dots, v_p) \oplus D\oplus W$ in $\PBC_\circ(V,U,W)$ is a preimage of $\{v_0,\dots,v_p\}$ under $\theta$, and $\theta$ is surjective. 

It remains to verify the third condition of  \autoref{defjoin}. Let $(v_0,C_0),\ldots, (v_p,C_p)$ be vertices of $\PBC_\circ(V,U,W)$ and let $\sigma_0, \dots,  \sigma_p$ be simplices  of $\PBC_\circ(V,U,W)$ such that for each $i=0, \ldots, p$, the vertex $(v_i,C_i)$ is a vertex of $\sigma_i$, and $\theta(\sigma_i)$ has vertices $\theta(v_0,C_0),\ldots, \theta(v_p,C_p)$, that is, vertices $\{v_0, \dots, v_p\}$.  We wish to show that the vertices $(v_0,C_0),\ldots, (v_p,C_p)$ span a simplex in $\PBC_\circ(V,U,W)$. 
 
 By  \autoref{alternatedescription}, it suffices to check that $v_i \in C_j$ for all $i\neq j$.  By assumption, $(v_i,D)$ for some $D\subseteq R^n$ and $(v_j, C_j)$ are vertices of $\sigma_j$. So a second application of  \autoref{alternatedescription} implies that $v_i \in C_j$ as required.

We now prove the simplicial complex $\PBC_\circ(V,U,W)$ is $ \left \lfloor \frac{\rk U -\rk W-s-2}{2}\right \rfloor$--connected. This follows by applying  \autoref{labelsystemconn} to the map $\theta$.  Since $\theta$ exihibits $\PBC_\circ(V,U,W)$ as a join complex over $\PB_\circ(U,W)$ and  $\PB_\circ(U,W)$ is wCM of dimension $(\rk U-\rk W-s)$ by \autoref{Vancor}, it only remains to check that the images of links of $p$--simplices are wCM of dimension $(\rk U -\rk W-s-p-2)$.

Let $\sigma = \{ (v_0, C_0), \dots, (v_p,C_p)\}$  be a $p$--simplex of $\PBC_\circ(V,U,W)$ and $C = \bigcap_{i=0}^p C_i$.  \autoref{linksPBC} \autoref{item:linkPBC} implies that
\[ \Lk_\circ (\sigma) = \PBC_\circ(C,U\cap C,W).\]
Note that
\[ \theta(\Lk_\circ(\sigma))=\PB_\circ(U\cap C,W) \subseteq \PB_\circ(U,W).\] Since $\rk(U \cap C) \ge \rk U-p-1$,  \autoref{Vancor} implies that $\theta(\Lk_\circ(\sigma))$ is wCM of dimension $(\rk U-\rk W-(p+1)-s)$, in particular wCM of dimension $(\rk U - \rk W - s - p - 2)$. 
%
%Combining  \autoref{linksPBC} and  this connectivity result  shows that $\PBC_\circ(V,U,W)$ is wCM of dimension $ \left \lfloor \frac{\rk U-\rk W-s}{2}\right \rfloor$. 
We conclude that $\PBC_\circ(V,U,W)$ is wCM of dimension $ \left \lfloor \frac{\rk U-\rk W-s}{2}\right \rfloor$. 
\end{proof}

 \autoref{CorPBCwCM} and  \autoref{complextosemi} (b) together establish  \autoref{PBCsemiconn}.

\subsection{Symplectic partial bases complexes} \label{SecSPB}

In this subsection, we consider the symplectic group analogues of the complexes from the previous section. Let $R$ be a PID. The high connectivity results discussed here will later be used to show that $\SI(R)$--modules satisfying a polynomial condition exhibit central stability.

\begin{definition} Let $W$ be a free finite-rank $R$--module equipped with a skew-symmetric bilinear form (not necessarily alternating or perfect). 
A \emph{symplectic partial basis} of $W$ is a set of pairs \[\{(v_0,w_0), \dots, (v_p,w_p)\}\subseteq W\times W\] such that $\{v_0,w_0, \dots, v_p,w_p\}$ is a partial basis of $W$ with
\[ \langle v_i,v_j \rangle = \langle w_i,  w_j\rangle = 0 \quad \text{and} \quad  \langle v_i, w_j\rangle = \delta_{ij}.\]
An \emph{ordered symplectic partial basis} is a symplectic partial basis with a choice of bijection between $\{(v_0,w_0), \dots, (v_p,w_p)\}$ and $\{0, \dots, p\}$.
%%%
\end{definition}

When $W$ is a symplectic $R$--module (that is, the form is alternating and perfect), then any symplectic partial basis may be extended to a symplectic basis for $W$. For general $W$, such a basis may not exist. %Following Friedrich \cite{Friedrich}, we make the following definition. 

\begin{definition}
Let $W$ be a free finite-rank $R$--module equipped with a skew-symmetric bilinear form. We call a submodule $H$ of $W$ \emph{symplectic} if the bilinear form restricts to a perfect alternating form on $H$. Define the \emph{Witt index} $g(W)$ of $W$ to be 
$$ g(W) =  \sup\left\{ {\textstyle \frac12} \, \rk(H)   \; \middle| \;  \begin{array}{ll} \text{there is a submodule $P \subseteq W$ and a symplectic} \\ \text{submodule $H \subseteq W$ such that } W \cong P \oplus H \end{array} \right\}.$$ 
.
\end{definition}

%\begin{definition}
%For $p \ge 0$, let $V$ be a free symplectic $R$--module. Let $\SPB_p(V)$ be the set of ordered symplectic partial bases of size $p+1$. There are maps $d_i\colon \SPB_p(V) \to %\SPB_{p-1}(V)$ given by forgetting the $i$th pair. With these maps, the sets assemble into a semisimplicial set $\SPB_\bullet(V)$. 
%\end{definition}

\begin{definition} Let $W$ be a free finite-rank $R$--module equipped with a skew-symmetric bilinear form. 
For $p \geq 0$, let $\SPB_p(V)$ be the set of ordered symplectic partial bases of size $p+1$. There are maps $d_i\colon \SPB_p(V) \to \SPB_{p-1}(V)$ given by forgetting the $i$th pair. With these maps, the sets assemble into a semisimplicial set $\SPB_\bullet(V)$. 
\end{definition} 

In the case that $W$ is a symplectic $R$--module, these semisimplicial sets (or the corresponding posets) have been studied, for example, by Charney \cite{Charney}, Panin \cite{Panin}, and Mirzaii--van der Kallen \cite{MvdK}. This generalization to arbitrary $W$ is due to Friedrich \cite{Friedrich}. In \cite{MvdK,Friedrich} more generally unitary groups are considered. In their notation, we are considering the case that the anti-involution is the identity, $\varepsilon=-1$, $\Lambda=R$, $\mu=0$, and $\lambda$  the symplectic form  is
\begin{align*} R^{2m} \times R^{2m} &\longrightarrow R, \\ \langle x,y\rangle &= \sum_{i=1}^m (x_{2i-1}y_{2i} -x_{2i}y_{2i-1}).\end{align*}
%Observe that
%\[ \SPB_p(V) \cap \Lk^{\MPB(V)}_\bullet (W) = \{ ((v_0,w_0), \dots, (v_p,w_p)) \in \SPB_p(V) \mid W \subseteq \spn(v_0,w_0, \dots, v_p,w_p)^\perp \}.\]
Mirzaii--van der Kallen \cite[Definition 6.3]{MvdK} defined a broad concept of the \emph{unitary stable rank} of a ring $R$, which specializes to apply to the study of symplectic, orthogonal, and unitary groups over $R$. 
%
%The following definition is the version of their definition relevant to symplectic forms.
%\begin{definition} \label{DefnUnitaryStableRank} A commutative ring $R$ has  \emph{(symplectic) unitary stable rank $s$} if $s$ is the smallest value of $m$ for which $R$ satisfies
%Condition (B$_m$) of \autoref{DefnStableRank}, and \emph{Condition (T$_{m+1}$)}. To state Condition (T$_m$) we define the group of  \emph{elementary symplectic matrices} $ES_{2m}(R)$ as follows.  Let $e_{i,j}(r)$ denote the $(2m \times 2m)$ matrix with $r \in R$ in position $(i,j)$ and 0 elsewhere. Then $ES_{2m}(R)$ is the group generated by the matrices of the form: 
%\begin{align*} I_{2m} + e_{2i-1, 2i}(r),  \; \; I_{2m}+e_{2i, 2i-1}(r), \; \; I_{2m} +e_{2i-1, 2j-1}(r) + e_{2j, 2i}(-r),  \\
% I_{2m} + e_{2i-1, 2j}(r) + e_{2j-1, 2i}(r), \; \;  I_{2m} + e_{2i, 2j-1}(r) + e_{2j, 2i-1}(-r) , && 1 \le i,j \le m, \; i \neq j. \end{align*}
 %  With this notation, Condition (T$_m$) is the statement that the group $ES_{2m}(R)$  acts transitively on the set of unimodular elements $v$ in $R^{2m}$. 
%\end{definition}
%
Results of Mirzaii--van der Kallen \cite[Remark 6.4]{MvdK} and Magurn--van der Kallen--Vaserstein \cite[Theorems 1.3 and 2.4]{MvdKV} together imply that semi-local rings (including fields) have symplectic unitary stable rank $1$, and PIDs have symplectic unitary stable rank at most $2$. Most relevant to this paper are the integers $R=\Z$, which have symplectic unitary stable rank $2$. 

%Having stated  \autoref{DefnUnitaryStableRank}, 
 Let $R$ be a PID of (symplectic) unitary stable rank $s$. Mirzaii and van der Kallen \cite[Theorem 7.4]{MvdK} proved that if $V \cong R^{2n}$ is a symplectic $R$--module, then  $\SPB_\bullet(V)$ is  $\left \lfloor \frac12\left({n-s-3}\right) \right \rfloor$--connected.  Analogously to the case of partial bases discussed in the previous subsection, Mirzaii and van der Kallen considered a simplicial complex which is homeomorphic but not isomorphic to the one considered here. Friedrich \cite[Theorem 3.4]{Friedrich} proved that for a general finite-rank free $R$--module $W$ with a skew-symmetric bilinear form of Witt index $g(W)$, then $\SPB_\bullet(W)$ is $\left \lfloor \frac12\left({g(W)-s-3}\right) \right \rfloor$--connected.

\section{Modules over stability categories}
\label{sec3}
In \autoref{SectionCategories}, we recall the definition of the categories $\VIC^H(R)$, $\SI(R)$, and $\FI$. Our approach to these topics uses the formalism of stability categories. We review the notions of central stability homology, degree of presentation, and polynomial degree for modules over these categories in \hyperref[centstabhom]{Sections \ref{centstabhom}},  \ref{SectionPolyDeg}, and \ref{SectionFinPolyDeg}.  The main result of \autoref{sec3} is that modules over $\VIC^H(R)$ and $\SI(R)$ with finite polynomial degree exhibit central stability, and we prove this result in \autoref{SectionMainLemma}. In \autoref{SectionSpecSequence},
we describe a spectral sequence introduced by Putman--Sam \cite{PutmanSamLinearGroups}, which we use in \autoref{secap}  to prove our representation stability results.

\subsection{Preliminaries} \label{SectionCategories}

In this subsection, we recall the definition of the category $\FI$ of Church--Ellenberg--Farb \cite{CEF} and the categories $\VIC^H(R)$ and $\SI(R)$ of Putman--Sam \cite[Section 1.2]{PutmanSamLinearGroups}. We will view these constructions as (equivalent to) stability categories (see \autoref{def:stabcatdef}), following the second author \cite[Section 3]{Pa2}. Stability categories are \emph{homogeneous categories} in the sense of Randal-Williams--Wahl \cite[Definition 1.3]{RWW} and \emph{weakly complemented categories} in the sense of Putman--Sam \cite[Section 1.3]{PutmanSamLinearGroups}. We will state their definition using a related concept, stability groupoids.

\begin{definition}\label{Def:stability groupoid}
Let $(\G,\oplus,0)$ be a monoidal groupoid whose monoid of objects is the natural numbers $\N$. The automorphism group of the object $n \in \N$ is denoted $G_n = \Aut^\G(n)$. Then $\G$ is called a \emph{stability groupoid} if it satisfies the following properties.
\begin{enumerate}
\item The monoidal structure
\[ \oplus \colon G_m\times G_n \inject G_{m+n}\]
is injective for all $m,n\in \N$.
\item The group $G_0$ is trivial.
\item $(G_{l+m}\times 1)\cap (1\times G_{m+n}) = 1\times G_{m} \times 1\subseteq G_{l+m+n}$ for all $l,m,n\in\mathbb N_0$.
\end{enumerate}
\end{definition}

\begin{definition}
A \emph{homomorphism} $ \mathcal G \to \mathcal H$ of stability groupoids is a strict monoidal functor sending $1$ to $1$.
\end{definition}

The construction used to define stability categories is originally due to Quillen but the notation and content of the following proposition is adopted from Randal-Williams--Wahl \cite[Section 1.1]{RWW}. 

\begin{proposition}[Randal-Williams--Wahl {\cite[Remark 1.4, Proposition 1.8]{RWW}}]\label{prop:symmetric monoidal}  
Let $\G$ be a stability groupoid. If $\G$ is braided, there is a monoidal category $U\G$ on the same objects, such that
\[ \Hom^{U\G}(m,n) \cong G_n/(1\times G_{n-m}).\]
as a $G_n$--set.

If a stability groupoid $\mathcal G$ is symmetric monoidal, then so is $U{\mathcal G}$.
\end{proposition}

\begin{definition}[Patzt {\cite[Definition 3.5]{Pa2}}]\label{def:stabcatdef}
Let $\G$ be a stability groupoid. If $\G$ is braided, then $U\G$ is called its \emph{stability category}.
\end{definition}

The following properties of stability categories are implicitly used throughout this section.

\begin{proposition}[Randal-Williams--Wahl {\cite[Proposition 1.7, 1.8(i)]{RWW}}]
Let $\G$ be a braided stability groupoid and $U\G$ its stability category, then:
\begin{enumerate}
\item $\G$ is the underlying groupoid of $U\G$.
\item The object $0\in U\G$ is initial. We will denote the unique map from $0$ to $n$ by $\iota_n$.
\end{enumerate}
\end{proposition}

\begin{example}\label{liststabcat} The following are examples of stability groupoids $\G$ and associated categories equivalent to $U\G$.
\begin{enumerate}    
\item \textbf{Symmetric groups:} Let $\mathfrak S$ be the symmetric stability groupoid of symmetric groups. Then $U\mathfrak S$ is equivalent to the category $\FI$  of finite sets and injection of Church--Ellenberg--Farb \cite{CEF}. \label{FIdef} 
\item \textbf{General linear groups:} Let $R$ be a ring. Let $\GL(R)$ be the symmetric stability groupoid of the general linear groups over $R$. Then $U\GL(R)$ is equivalent to the category $\VIC(R)$ of Djament \cite{DjamentUnitaires} and Putman--Sam \cite{PutmanSamLinearGroups}. \label{defVIC} The objects of $\VIC(R)$ are finite rank free $R$--modules and the morphisms from $V$ to $W$ are given by the set of pairs $(f,C)$ of an injective homomorphism $f\colon V \to W$ and a free submodule $C \subseteq W$ such that $\im(f) \oplus C=W$ and $\rk(C) + \rk(V)=\rk(W)$.  The composition law is defined by \[(f,C) \circ (g,D) = (f \circ g,C \oplus f(D)).\] A \emph{$\VIC(R)$--module} is a functor from $\VIC(R)$ to the category of abelian groups.

\item \textbf{General linear groups with restricted determinant:} Let $R$ be a commutative ring and $H$ a subgroup of the group of units $R^\times$. Let $\GL^H(R)$ denote the stability groupoid of the subgroups
\[ \GL^H_n(R) = \{ A \in \GL_n(R) \mid \det A \in H\}.\]
$\GL^H(R)$ is symmetric monoidal if $-1\in H$. Note that determinant is well defined since the ring is commutative. Its stability category $U\GL^H(R)$ is equivalent to the category $\VIC^H(R)$, defined by Putman--Sam \cite[Section 1.2]{PutmanSamLinearGroups}. Its objects are finite-rank free $R$--modules $V$ with a choice of generator of $a \in \left( \bigwedge^{\rk(V)} V \right) /H$. The morphisms $\Hom^{\VIC^H(R)}((V,a),(W,b))$ are split injective linear maps $f:V \m W$ and a choice of complement of the image of $f$ with the added requirement that $f_*(a)=b$ if $f$ is an isomorphism. Here $f_*: \left( \bigwedge^{\rk(V)} V \right) /H \m \left( \bigwedge^{\rk(W)} W \right) /H$ is the induced map. In particular, $$\Hom^{\VIC^H(R)}((V,a),(W,b)) \cong  \Hom^{\VIC(R)}(V,W), \qquad \text{ for $\rk(V) \neq \rk(W)$, and}$$ $$\Hom^{\VIC^H(R)}((V,a),(W,b)) \cong \GL_n^H(R) \qquad \text{ if $V \cong W \cong R^n$.}$$

% are again injective complemented linear maps $(f,C)$, and we assign to $C$ the (unique) $H$--orientation that makes the $H$--orientations on $(\im f \oplus C)$ and $W$ agree. 
% In particular, linear maps in  $\Hom^{\VIC^H(R)}(V,V)$ must have determinants in $H$.

\item \textbf{Symplectic groups:} Let $R$ be a ring. Let $\Sp(R)$ be the symmetric stability groupoid of the symplectic groups over $R$. Then $U\Sp(R)$ is equivalent to the category $\SI(R)$ of free finite-rank symplectic $R$--modules and isometric embeddings.  Details are given in  Putman--Sam \cite[Section 1.2]{PutmanSamLinearGroups}.
\item \textbf{Automorphisms of free groups:} Let $\AutF$ be the symmetric stability groupoid of the automorphism groups of free groups of finite rank. Then $U\AutF$ is equivalent to the category of finite-rank free groups and monomorphisms together with a choice of free complement, that is, 
\[ \Hom(F,G) = \left\{ (f,C) \middle| \begin{array}{l} f\colon F \inject G \text{ an injective group homomorphism,} \\ \text{$C \subseteq G$ a free subgroup with }  \; (\im f )* C = G \end{array} \right\}.\]
See Randal-Williams--Wahl \cite[Section 5.2.1]{RWW}.
\item \textbf{Mapping class groups of compact oriented surfaces with one boundary component:} Let $\Mod$ be the braided stability groupoid of mapping class groups of compact oriented surfaces with one boundary component.  Its monoidal structure is induced by boundary connect sum. See Randal-Williams--Wahl \cite[Section 5.6.3]{RWW}. 
\end{enumerate} 
\end{example}

\begin{definition} We call functors from a category $\C$ to the category of abelian groups \emph{$\C$--modules} and denote the category of $\C$--modules by $\Mod_\C$. If $\C$ is a stability category or stability groupoid and $A$ is a $\C$--module, we let $A_n$ denote the value of $A$ on the object $n\in \N$.

\end{definition}

If a category $\mathcal{C}$ is equivalent to the category $U\G$, then the categories of $\mathcal{C}$--modules and $U\G$--modules are equivalent. We therefore use the terms of $U\mathfrak S$--module and $\FI$--module roughly interchangeably, and similarly for other items in \autoref{liststabcat}.

%We now describe a class of particularly simple $U\G$--modules. 

%Let in the following $\G$ be a braided stability groupoid with the automorphism groups $\{G_n\}_{n\in\N}$ and $U\G$ its stability category.

%Note that a $\G$--module (a functor from $\G$ to abelian groups) is precisely the data of a $\Z[\G_n]$--module for every $n$. Let $\Mod_{U\G}$ denote the category of $U\G$--modules, $\Mod_{\G}$ denote the category of $\G$--modules and $F:\Mod_{U\G} \m \Mod_{\G}$ denote the forgetful functor. 

If $\G$ is a stability groupoid, then the data of $\G$--module is precisely the data of a $\Z[G_n]$--module for every $n$. For $m$ fixed, we will often view a $\Z[G_m]$--module $W$ as a $\G$--module by putting $W$ in degree $m$ and the module $0$ in all other degrees. We now define free $U\G$--modules. 

%Given a $\G$--module $W$, let $W_n$ denote the associated $\Z[\G_n]$--module. 

\begin{definition}
Let $\G$ be a stability groupoid and let $M:\Mod_{\G} \m \Mod_{U\G}$ be the left adjoint to the forgetful functor $\Mod_{U\G} \m \Mod_{\G}$. We say that $M(W)$ is the the \emph{free} $U\G$--module on $W$. Concretely, given a $\Z[G_m]$--module $W$, the $\Z[G_n]$--module $M(W)_n$ is given by the formula
\[M(W)_n \cong  \Z[\Hom^{U\G}(m,n)] \otimes_{\Z[G_m]} W . \] 
%\text{where $\Z$ is the trivial $\Z[G_{n-m}]$--module.} \end{array} \right. \] 
For a general $\G$--module  $W$, \[M(W) = \bigoplus_{n=0}^{\infty} M(W_n).\]  We abbreviate $M(\Z[G_m])$ by $M(m)$. \label{deffree}
\end{definition}

\subsection{Central stability homology and resolutions} \label{centstabhom}
We begin by defining degree of generation and presentation and discuss how these concepts relate to central stability degree. We then review central stability homology and how it relates to the degrees of higher syzygies of $U\G$--modules.

\begin{definition}\label{gen} Let $\G$ be a stability groupoid. We say a $\G$--module $W$ has \emph{degree} $\le d$ if $W_n \cong 0$ for $n>d$. A $U\G$--module $A$ is \emph{generated in degrees $\le d$} if there is an exact sequence of $U\G$--modules 
\[  M(W^0) \m A \m 0 \qquad \qquad \] with $W^0$ of degree $\le d$. A $U\G$--module $A$ is \emph{related in degrees $\le d$} if there is an exact sequence of $U\G$--modules 
\[  M(W^1) \m M(W^0) \m A \m 0 \qquad \qquad \] with $W^1$ of degree $\le d$. We say $A$ is \emph{presented in degrees $\le d$} if it is generated and related in degrees $\le d$. 
%If for each $n\in \N$ the abelian group $A_n$ is finitely generated, then \emph{generation in finite rank} is  equivalent to the concept of \emph{finite generation} that is central to the work of Church--Ellenberg--Farb \cite{CEF}. 
\end{definition}

\begin{definition}\label{def:centralstability}
Let $\G$ be a braided stability groupoid with braiding $b$. Let $A$ be a $U\G$--module, and let $\phi_{n-2} \colon A_{n-2} \to A_{n-1}$ denote the map induced by the morphism
\[ \id_{n-2}\oplus \iota_1\colon (n-2)\oplus 0 \longrightarrow (n-2)\oplus 1\]
that is the sum of the identity map of $n-2$ and the unique morphism $0 \to 1$. We consider two maps
\[ \Ind^{G_n}_{G_{n-2}\times 1} A_{n-2} \rightrightarrows \Ind^{G_n}_{G_{n-1}\times 1} A_{n-1}.\]
The first is induced by $\phi_{n-2}$ and the second is the map induced by $\phi_{n-2}$ postcomposed with $\id_{n-2} \oplus (b_{1,1})$.

\noindent The $U\G$--module $A$ has \emph{central stability degree $\le d$} if the map
\[ \coeq( \Ind^{G_n}_{G_{n-2}} A_{n-2} \rightrightarrows \Ind^{G_n}_{G_{n-1}} A_{n-1}) \longrightarrow A_n\]
induced by $\phi_{n-1}$ is an isomorphism for all $n>d$. %We say $A$ is \emph{centrally stable} if it has finite central stability degree.
\end{definition}

%The following remark gives a word of warning towards the definition of central stability. Despite this warning the reader may be assured that for all categories considered in this paper, central stability is equivalent to presentability in finite degree. Thus all definitions found in the literature coincide in these cases. 

%The statement concerning $\FI$--modules can be deduced from Church--Ellenberg \cite[Proposition 4.2]{CE} and the statement involving other categories is due to Patzt \cite[Corollary 6.4]{Pa2}. 

\begin{remark} \label{RemarkCSWarning}
The definition of central stability used by Putman \cite[see paragraph ``Central stability, definition"]{Pu}, the one used by Putman--Sam \cite[see paragraph ``Asymptotic structure"]{PutmanSamLinearGroups} and the one we use are all slightly different. These concepts turn out to coincide for the symmetric groups \cite[Remark 1.2 \& Theorem F]{PutmanSamLinearGroups}. The work of Djament \cite[Proposition 2.14]{DjamentFoncteurs} implies that Putman--Sam's definition agrees with the definition of degree of presentation used in this paper; see also Gan--Li \cite[Theorem 3.2]{GanLiCentralStability}. It is unclear how to generalize Putman's definitions to braided stability categories, but there is a straightforward generalization to symmetric stability categories, which is equivalent to our definition. See \cite[Proposition 6.1 and 6.2]{Pa2} for a comparison of these different forms of stability.
\end{remark}

We now compare the notion of central stability degree to the notion of presentation degree. The following theorem will allow us to translate results on presentation degree quoted from other papers into statements about central stability. This theorem also allows us to translate the main theorems of this paper into results about presentation degree.

%\begin{theorem}[Patzt {\cite[Corollary \todo{???}]{Pa2}}]

%\label{thm:centralpresent}   For $\G=\mathfrak S$, $\GL$, or $\SI$, \[\text{let }(k,a)=\begin{cases}
%(1,1) & \text{if } \G=\mathfrak S \\
%(2,s+1), & \text{if } \G=\GL^H(R) \text{ with } R \text{ a ring of stable } s \\
%(2,s+2), & \text{if } \G=\SI(R) \text{ with } R \text{ a ring of unitary stable } s 
%\end{cases}\] \noindent A $U\G$--module $A$ be a with $g,d \in \N$ and take $d  \ge 2 \max(k,a)$.  Then $A$ is presented in degrees $\le d$ if and only if $A$ has central stability degree $\le d$.

%\end{theorem}

%\begin{theorem}[Church--Ellenberg {\cite[Proposition 4.2]{CE}}, Patzt {\cite[Corollary \todo{???}]{Pa2}}] 

%\label{thm:centralpresent}   For $\C=\FI$, $\VIC^H(R)$, or $\SI(R)$, let $d \in \N$ satisfy: \todo{Correct this statement.}
%\begin{itemize}
%\item $d \ge 0$ if $\C=\FI$
%\item $d \ge 2s+2$  if $ \C=\VIC^H(R)$ with $ R$  a ring of stable rank $ s $
%\item $d \ge 2s+4$  if $ \C=\SI(R)$ with $ R$  a ring of (symplectic) stable unitary rank stable $ s $.
%\end{itemize}
% \noindent Then a $\C$--module $A$ is presented in degrees $\le d$ if and only if $A$ has central stability degree $\le d$.
%\end{theorem}

\begin{theorem} \label{thm:centralpresent} Let  $A$ be a module over $\FI$, a module over $\VIC^H(R)$ for a ring $R$ of finite stable rank and a subgroup $H\le R^\times$ containing $-1$, or a module over $\SI(R)$ for a ring $R$ of finite (symplectic) unitary stable rank. Then $A$ is presented in finite degrees if and only if it has finite central stability degree. Specifically, we have the following bounds. 
\begin{description}
\item[$\FI$ \textnormal{(Patzt {\cite[Corollary 6.4(a)]{Pa2}}; see also Church--Ellenberg {\cite[Proposition 4.2]{CE}})}]  For $d \ge 0$, let  $A$  be an $\FI$--module generated in degrees $\le d$. Then $A$ is presented in degrees $\leq (d+1)$ if and only if $A$ has central stability degree $\le (d+1)$. \label{FIcentpresnet} 
%[CE, Prop 4.2] only proves one direction, but the other direction is also true, because \widetilde H_0 of free FI-modules is always zero (in all degrees)

\item[$\VIC^H(R)$ \textnormal{(Patzt {\cite[Corollary 6.4(b)]{Pa2}})}]  Suppose $R$ is a ring of stable rank $s$, $-1\in H\le R^\times$, and let $A$ be a $\VIC^H(R)$--module generated in degrees $\le g$. If $A$ has central stability degree $\le d$, then it is presented in degrees $\le \max\big(d, g+s+1\big)$. If $A$ is related in degrees $\le d$, then it has central stability degree $\le \max\big(d,g+s+1\big)$.

\item[$\SI(R)$ \textnormal{(Patzt {\cite[Corollary 6.4(c)]{Pa2}})}] Suppose $R$ is a ring of (symplectic) unitary stable rank $s$ and let $A$ be a $\SI(R)$--module generated in degrees $\le g$. If $A$ has central stability degree $\le d$, then it is presented in degrees $\le \max\big(d, g+s+2\big)$. If $A$ is related in degrees $\le d$, then it has central stability degree $\le \max\big(d, g+s+2\big)$.

\end{description}
\end{theorem}

We note that for $H \neq R^\times$, we implicitly assume $R$ is commutative when considering $\VIC^H(R)$ as we have not defined $\VIC^H(R)$ for noncommutative rings, only $\VIC(R)$. 

Observe that the central stability degree of a $U\G$--module is always at least its degree of generation. Thus the central stability degree alone gives a bound on degree of presentation, though this bound may be improved with the data of both central stability degree and the degree of generation. 

In fact, the following lemma states an easy equivalence that has been pointed out in many papers, for example \cite[Proposistion 5.4]{Pa2}.

\begin{lemma}\label{lem:generation}
A $U\G$--module $A$ is generated in degrees $\le d$ if and only if 
\[ \Ind_{G_{n-1} \times 1}^{G_n} A_{n-1} \longrightarrow A_n\]
is surjective for all $n>d$.
\end{lemma}

%To state the relationship between presentation degree and central stability degree, we need the following definition. 

Given a $U\G$--module $A$, central stability can be rephrased as acyclicity of a chain complex 
\[\Ind_{G_{n-2}}^{G_n} A_{n-2} \longrightarrow  \Ind_{G_{n-1}}^{G_n} A_{n-1} \longrightarrow  A_n \longrightarrow  0.\] This chain complex is the tail of a longer chain complex:
%the \emph{central stability chain complex} of $A$, which we describe in \autoref{DefnChainComplex}. We call its homology \emph{central stability homology}. \autoref{thm:centralpresent} is a special case of \autoref{thm:res of finite type} (quoting \cite[Theorem 5.7]{Pa2}) which states that vanishing of central stability homology controls the higher syzygies of $A$. 

%The second author \cite{Pa2} defines \emph{central stability homology}  for $\C$--modules when $\C$ is a monoidal category whose unit object is initial. For present purposes we specialize to the case when $\C$ is a stability category $U\G$, and define an augmented semisimplicial $U\G$--module $\widetilde C_\bullet(A)$ and associated chain complex  $\widetilde C_*(A)$  for any $U\G$--module $A$. 

\begin{definition}[Compare {\cite[Definition 2.4, Proposition 4.3]{Pa2}}] \label{DefnChainComplex} Let $A$ be a $U\G$--module, and let $\phi_n\colon A_n \to A_{n+1}$ denote the map induced by the morphism 
\[ \id_n \oplus \iota_1 \colon n \oplus 0 \longrightarrow n \oplus 1.\] 
Define $\widetilde C^\G_\bullet(A)$ to be the augmented semisimplicial $U\G$--module with $p$--simplices given by
\[ \widetilde C^\G_p(A)_n = \Z[G_n] \otimes_{\Z[G_{n-(p+1)} \times 1]} A_{n-(p+1)}\] and the $i$th face map by
\begin{align*}
 d_i \colon \Z[G_n] \otimes_{\Z[G_{n-(p+1)}\times 1]} A_{n-(p+1)} &\longrightarrow \Z[G_n] \otimes_{\Z[G_{n-p}\times 1]} A_{n-p}\\
 g \otimes a &\longmapsto g\cdot(\id_{n-p}\oplus h_i) \otimes \phi_{n-p-1}(a)
 \end{align*}
where the coset $h_iG_1$ corresponds to the morphism
\[ \id_i \oplus \iota_1 \oplus \id_{p-i} \colon i \oplus 0 \oplus (p-i) \longrightarrow i \oplus 1 \oplus (p-i) \quad  \text{ in } \Hom^{U\G}(p,p+1)\cong G_{p+1}/G_{(p+1)-p}.\]

Let $\widetilde C^\G_*(A)$ denote the chain complex associated to $\widetilde C^\G_\bullet(A)$; we will refer to it as the \emph{central stability complex} of $A$. Let $\widetilde H^\G_*(A)$ denote the homology of this chain complex; we will call it the \emph{central stability homology} of $A$.

\end{definition}

%For an $\FI$--module $A$, the chain complex $\widetilde C_*(A)$ is closely related to Putman's chain complex $\mathrm{IA}_{*+1}(A_{n-*-1})$ \cite[Section 4]{Pu}, and the complexes computing \emph{$\FI$--homology} in work of Church, Ellenberg, Farb, Nagpal \cite{CEF, CEFN, CE} and Gan \cite{GanLES}. In these complexes, the $p$-chains are the groups  $\widetilde C_p(A) \otimes_{\Z[\mathfrak{S}_{p+1}]} \Z$ where $\Z$ is the sign representation of  $\mathfrak{S}_{p+1}$.  

The chain complex $\widetilde C^\G_*(A)$ and groups $H^\G_i(A)$ have an $\N$-grading  coming from the fact that the set of objects of a stability groupoid is $\N$. We denote the piece in grading $n$ by $\widetilde C^\G_*(A)_n$ or $\widetilde H^\G_i(A)_n$ respectively. We will  drop the superscript $\G$ when the groupoid is clear from context.

\begin{definition} \label{highercentstab}
We say that a $U\G$--module $A$ has \emph{higher central stability} if $\widetilde H_i(A)_n \cong 0$ for $n$ sufficiently large compared to $i$. 
\end{definition}

\begin{remark}
Note that in particular $\widetilde H_{-1}(A)_n=0$ and $\widetilde H_{0}(A)_n=0$ is equivalent to the surjectivity and the injectivity, respectively, of the last map in \autoref{def:centralstability}.
\end{remark}

The central stability complex has appeared in the literature in varying degrees of generality. In the $\FI$ case, it is closely related to a complex introduced in \cite[Section 4]{Pu} whose homology is called \emph{$\FI$-homology} by Church-Ellenberg \cite{CE}. $\FI$--homology has appeared in work of Church, Ellenberg, Farb, and Nagpal \cite{CEF, CEFN, CE} and Gan and Li \cite{ GanLiCentralStability, GanLiRemark, GanLES, LiHomInvariants}. For the category $\FI$, the complex $\widetilde C_*(A)$ itself was denoted by $B_{*+1}(\mathcal A)$ by Church--Ellenberg--Farb--Nagpal \cite[Definition 2.16]{CEFN}, by $\widetilde{C}^{F_{\mathcal A}}_{*+1}$ by Church--Ellenberg \cite[Section 5.1]{CE}, and by $\mathrm{Inj}_*(A)$ by the first and third author \cite[Section 2.2]{MillerWilson}. Putman--Sam \cite[Section 3]{PutmanSamLinearGroups} defined the complex for modules $A$ over general ``cyclically generated''  \emph{complemented categories}  and used the notation $\Sigma_{*+1} (A)$. In this paper, we adopt notation used by the second author \cite{Pa2}.
% Miller--Wilson \cite[Definitions 2.19 and 2.46]{MillerWilson} called this complex $\mathrm{Inj}_*(A)$ for an $\FI$--module $A$, and $\mathrm{Inj}^2_*(A)$ when $A$ is a module over the \emph{twisted skew-commutative algebra} $\bigwedge \mathrm{Sym}^2 \Z$. 
%\todo{Cite Gan--Li \cite{GanLiCentralStability} better}

When $A=M(0)$, the central stability homology is the reduced homology of the following augmented semisimplicial set that was previously studied by Randal-Williams--Wahl \cite[Definition 2.2, Theorem 3.1, 4.20]{RWW}.

\begin{definition}
Let $K_\bullet \G$ be the semisimplicial $U\G$--set with set of $p$--simplices given by 
\[ (K_p\G)_n = \Hom^{U\G}(p+1,n) \cong G_n/G_{n-(p+1)} \] and the $i$th face map is induced by precomposition with the morphism
\[ \id_i \oplus \iota_1 \oplus \id_{p-i} \colon i \oplus 0 \oplus (p-i) \longrightarrow i \oplus 1 \oplus (p-i).\]
\end{definition}

By definition, the reduced chain complex  $\widetilde C_*((K_\bullet \G)_n)$ agrees with the central stability complex of $M(0)$. Here the  term reduced chain complex of a semisimplicial set means the chain complex formed by first augmenting the semisimplicial set by inserting a point in degree $-1$, then composing with the free abelian group functor to obtain an augmented semisimplicial abelian group, and then forming a chain complex by taking the alternating sum of face maps as the differential.

%\label{identificationRemark}
The semisimplicial sets $(K_\bullet \G)_n$ were first introduced by Randal-Williams--Wahl \cite{RWW} and are examples of a larger class of semisimplicial sets which they denote by $W_n(A,X)$. They describe these semisimplicial sets explicitly in many examples, often showing they are isomorphic to complexes previously considered in the literature. In particular, the  semisimplicial set $(K_\bullet \GL(R))_n$ is isomorphic to $\PBC_\bullet(R^n)$  \cite[Section 5.3]{RWW},  $(K_\bullet \Sp(R))_n$ is isomorphic to $\SPB_\bullet(R^n)$ \cite[Section 5.4]{RWW}, and $(K_\bullet \mathfrak S)_n$ is isomorphic to the \emph{complex of injective words} \cite[Section 5.1]{RWW}, a complex introduced by Farmer \cite{Fa}.  Randal-Williams--Wahl \cite{RWW} proved that high connectivity of $(K_\bullet \G)_n$ implies homological stability for the groups $G_n$. The name \textbf{H3} is related to their \cite[Definition 2.2]{RWW}.

\begin{definition}
Let $a,k\in \N$. Define the following condition for a stability category $U\G$. 
\begin{description}
\item[H3($k$,$a$)] $\widetilde H_i(M(0))_n = 0$ for all $i\ge-1$ and all $n> k\cdot i + a$.
\end{description}
\end{definition}

\begin{remark}\label{rem:H_{-1}}
For all stability categories $U\G$, $\widetilde H_{-1}(M(0))_n=0$ for all $n>0$.
\end{remark}

In the following proposition, we compile information from the literature about this condition for the stability categories appearing in  \autoref{liststabcat}.

\begin{proposition}\label{H3} These stability categories satisfy the condition \con{H3}.
\begin{enumerate}
\item $U\mathfrak S$ satisfies \con{H3($1$, $1$)}.
\item $U\AutF$ and $U\Mod$ satisfy \con{H3($2$, $2$)}. 
\item If $R$ is a ring with  stable rank $s$, then $U\GL(R)$ and satisfies \con{H3($2$, $s+1$)}. If $R$ is a commutative ring with  stable rank $s$, then $U\GL^H(R)$ and satisfies \con{H3($2$, $s+1$)}. 
\item If $R$ is a ring with (symplectic) unitary stable rank $s$, then $U\Sp(R)$ satisfies \con{H3($2$, $s+2$)}. 
\end{enumerate}
\end{proposition}

\begin{proof}
Since $\widetilde H_i((K_\bullet \G)_n) \cong \widetilde H_i(M(0))_n$, it suffices to show $(K_\bullet \G)_n$ is highly connected. With the exception of $U\GL^H(R)$, these connectivity results appear in Randal-Williams--Wahl \cite{RWW}. 
\begin{enumerate} 
\item Randal-Williams--Wahl \cite[Section 5.1]{RWW} proved that $(K_\bullet \mathfrak S)_n$ is isomorphic to complex of injective words which which was shown to be $(n-2)$--connected by Farmer \cite[Theorem 5]{Fa}.

\item Hatcher--Vogtmann \cite[Proposition 6.4]{HatcherVogtmannCerf} showed a space closely related to $K_\bullet U\AutF$ is highly connected. Similarly, Hatcher--Vogtmann \cite[Main Theorem]{HatcherVogtmannTethers} proved that a space closely related to $(K_\bullet U\Mod)_n$ is highly connected. Using these results,  Randal-Williams--Wahl \cite[Proposition 5.3 and Lemma 5.25]{RWW} proved that $(K_\bullet \AutF)_n$ and $(K_\bullet \Mod)_n$ are $ \left \lfloor \frac{1}{2}(n-3) \right \rfloor$--connected. 
\item Randal-Williams--Wahl \cite[Lemma 5.10]{RWW} proved  that $(K_\bullet \GL(R))_n$ is $ \left \lfloor \frac{1}{2}(n-s-2)\right \rfloor$--connected. This builds on the work of van der Kallen \cite[Theorem 2.6]{V}. Also see Charney \cite[Theorem 3.5]{Ch1}.
\item From the description of $(K_\bullet \Sp(R))_n$ in \cite[Section 5.4]{RWW}, it follows that it is isomorphic to the space appearing in Mirzaii--van der Kallen \cite[Theorem 7.4]{MvdK} and hence is %$ \displaystyle \left \lfloor \frac{ {\textstyle \frac12} \rk V-s-3}{2}\right \rfloor$--connected.
$  \left \lfloor \frac{1}{2}(n-s-3)\right \rfloor$--connected.
\end{enumerate}

%Also see Randal-Williams--Wahl \cite[Section 5]{RWW} for a discussion of these complexes from a cat. 

%The results for all categories but $U\GL^H$ appear in Randal-Williams--Wahl \cite[Section 5]{RWW}. 

We now consider the case $\G=\GL^H(R)$. Because
\[ \Hom^{U\GL^H(R)}(m,n) = \Hom^{U\GL(R)}(m,n) \qquad  \text{for $m<n$,}\]
it follows that 
\[ \widetilde C_p^{\GL^H(R)}(M(0))_n = \widetilde C_p^{\GL(R)}(M(0))_n \qquad \text{for $p\le n-2$.} \]
 We conclude that the map
\[ \widetilde C_*^{\GL^H(R)}(M(0))_n \inject \widetilde C_*^{\GL(R)}(M(0))_n \]
induces an isomorphism on homology groups for $*\le(n-3)$. Thus for $i\ge 0$,
\[  \widetilde H_i^{\GL^H(R)}(M(0))_n = \widetilde H_i^{\GL(R)}(M(0))_n =0 \qquad \text{when $n>2i+s+1$}\]
since $U\GL(R)$ satisfies \con{H3($2$, $s+1$)} and necessarily $(n-3) \ge i$ in this range. Moreover, as in \autoref{rem:H_{-1}},
\[  \widetilde H_{-1}^{\GL^H(R)}(M(0))_n = 0 \qquad \text{for $n>0$}\]
and we conclude \con{H3($2$, $s+1$)} for $U\GL^H$.
\end{proof}

The following theorem generalizes \autoref{thm:centralpresent} to general stability categories and higher central stability homology groups.

\begin{theorem}[Patzt {\cite[Theorem 5.7]{Pa2}}]\label{thm:res of finite type}
Assume \con{H3($k$,$a$)}. Let $A$ be a $U\G$--module and $\{d_n\}_{n\in\N}$ a sequence of integers with $d_{i+1}-d_i \ge \max(k,a)$,
then the following statements are equivalent.
\begin{enumerate}
\item \label{ItemFreeResolution} There is a resolution
\[ \dots \to M(W^1) \to M(W^0) \to A \to 0\]
with $W^i$ a $\G$--module of degree $\le d_i$.
\item \label{ItemVanishingHomology} The homology groups
\[ \widetilde H_i(A)_n \cong 0 \qquad \qquad  \text{for all $i\ge -1$ and all $n > d_{i+1}$.} \]
\end{enumerate}
\end{theorem}

\subsection{Polynomial degree} \label{SectionPolyDeg}

%We achieve a main goal of this paper -- establishing central stability of the second homology groups of Torelli groups -- by proving vanishing of their central stability homology.

To prove that the second homology of Torelli subgroups are centrally stable, will first prove that the first homology group of these Torelli groups exhibit higher central stability. 
%The behavior of central stability homology groups is not yet well understood in general. For example, given a $\VIC(\Z)$ or $\SI(\Z)$--module $A$ with finite central stability degree, it is currently unknown whether $A$ has higher central stability. Indeed, this result would imply central stability for the homology of the corresponding Torelli groups in every homological degree. We note that over a field of characteristic zero, it is true that central stability for $\VIC(\Z/p\Z)$ or $\SI(\Z/p\Z)$--modules implies higher central stability \cite[Theorems C and D]{MillerWilson2}.  
To do this, we establish higher central stability for \emph{polynomial} $\VIC(\Z)$-- and $\SI(\Z)$--modules. We can apply these results to the Torelli groups because their first homology groups are known to be polynomial. There are several notions of polynomial functors---as mentioned in the introduction. We use the same definition as in the arXiv version \cite[Definition 4.10]{RWWarxiv}; in the published version the definition was slightly changed.

%The main technical result of the paper is the following: 
 %for a stability category $U\G$ satisfying the Condition \con{H4} (\autoref{H4}), in each degree $i$, the central stability homology $\widetilde H_i(A)_n$ of a $U\G$--module $A$ is eventually zero if $A$ has finite polynomial degree. This criterion is often easy to check in practice, and we use it to prove our stability results for the Torelli groups and congruence subgroups. 

\begin{definition} Let $\G$ be a stability groupoid. Define the endofunctor 
\begin{align*} S& \colon U\G\to U\G \\ \text{via the formula \, \,} S&= 1\oplus -. \end{align*} 
We will consider the natural transformation $\id \to S$ given by \[ \iota_1 \oplus \id_n : 0 \oplus n \longrightarrow 1\oplus n.\] 
By abuse of notation, we denote the endofunctor of $U\G$--modules given by precomposition by $S$ also by $S$. Again, there is an induced natural transformation $\id \to S$ defined by precomposition with the above natural transformation. 
\end{definition}
Concretely, if $\G$ is braided and $A$ is a module over its stability category $U\G$, then there are isomorphisms of $G_n$--representations 
\[(SA)_n \cong \mathrm{Res}^{G_{n+1}}_{1 \times G_n} A_{n+1} \cong  \mathrm{Res}^{G_{n+1}}_{ G_n \times 1} A_{n+1}  .\] This last isomorphism uses the braiding. 

\begin{definition} Given a $U\G$--module $A$,  define $U\G$--modules \[\qquad \quad  \ker A := \ker(A \to SA) \qquad \text{and} \]  \[\coker A := \coker(A\to SA).\] 
\end{definition}

\begin{definition} \label{defpoly}
We say that $A$ has \emph{polynomial degree $-\infty$ in ranks $>d$} if $A_n = 0$ for all $n>d$. We say $A$ has \emph{polynomial degree $\le 0$ in ranks $>d$} if $\ker A_n = 0$ for all $n>d$  and $\coker A$ has polynomial degree $-\infty$ in ranks $>d$. For $r\ge1$, we say $A$ has \emph{polynomial degree $\le r$ in ranks $>d$} if $\ker A_n = 0$ for all $n>d$  and $\coker A$ has polynomial degree $\le r-1$ in ranks $>d$. 

We say $A$ has \emph{polynomial degree $\le r$} if it has polynomial degree $\le r$ in all ranks $>-1$.
\end{definition}

%\begin{remark} The indexing convention used in \autoref{defpoly} is consistent with work of the second author \cite{Pa2} but differs from the convention used by Randal-Williams--Wahl \cite{RWW}. A $U\G$--module $A$ that has polynomial degree $\le r$ in ranks $>d$ in the sense of \autoref{defpoly} is a \emph{coefficient system of  degree $\le r$ at $(d+r)$} in the sense of Randal-Williams--Wahl \cite[Definition 4.10]{RWW}. \todo{This version of RWW is not on the ArXiv yet} \end{remark}

Polynomial modules have higher central stability when the category satisfies the following condition.

\begin{definition}\label{H4}
Let $b,\ell\ge 1$ be natural numbers. Define the following condition on a stability category $U\G$: 
\begin{description}
\item[H4($\ell$,$b$)] $\widetilde H_i(\coker M(m))_n = 0$ for all $m \ge 0$, all $i\ge-1$, and all $n> \ell\cdot (i+m) + b$.
\end{description}
\end{definition}

\begin{theorem}[Patzt {\cite[Corollary 7.9]{Pa2}}]\label{thm:findeg homology}\label{cor:findeg homology}
Let $a,b,k,\ell\ge 1$. Let $\G$ be a braided stability groupoid, and assume its stability category $U\G$ satisfies \con{H3($k$,$a$)} and \con{H4($\ell$,$b$)} with $b\ge\max(k,a)$. 
If $A$ is a $U\G$--module of polynomial degree $\le 0$ in ranks $>d$ for some $d \geq -1$, then \[\widetilde H_i(A)_n = 0 \qquad \text{ for all \quad $i\ge -1$ \quad and \quad  $n> \max(d+i+2, ki+a)$.} \]
If $A$ is a $U\G$--module of polynomial degree $\le r$ in ranks $>d$ for some $r \ge 1$ and $d \geq -1$, then \[\widetilde H_i(A)_n = 0 \qquad \text{ for all \quad $i\ge -1$ \quad and \quad  $n> \ell^{i+1}(d+r) + (\ell^i + \ell^{i-1} + \dots + \ell + 1)b+1$.} \]
\end{theorem}

Note that for $i=-1$, the above sum is an empty sum and hence zero.

\begin{proof}[Proof of \autoref{thm:findeg homology}]
\cite[Corollary 7.9]{Pa2} only considers the case $r\ge 1$. But the cases $r\le 0$ are easy to see. First note that if $A_n =0$ for all $n>d$, also $\widetilde C_i(A)_n \cong \Ind_{G_{n-i-1}}^{G_n} A_{n-i-1} =0$ for all $n>i+1+d$ and thus $\widetilde H_i(A)_n = 0$ for all $n>i+1+d$.

Now let us assume $A$ has polynomial degree $\le 0$ in ranks $>d$. From the definition of polynomial degree, it follows that for $n > d$, the maps $A_n \to (SA)_n$ are isomorphisms. Thus there is a short exact sequence of $U\G$--modules
\[ 0\longrightarrow B \longrightarrow A \longrightarrow C \longrightarrow 0\]
 such that $B_n =0$ for all $n\le d$, $B_n \to (SB)_{n}$ is an isomorphism for all $n>d$, and $C_n = 0$ for all $n>d$. Here $B$ is obtained from $A$ by replacing the groups in ranks $\leq d$ with $0$ and $C$ is the cokernel of the natural map $B \m A$.

We will show moreover that, for $n >d$, the group $G_n$ acts trivially on $A_n$. Fix $n>d$. %If $n=0$, then $G_0=1$. So suppose $n>0$, and consider 
Consider the map $\rho: A_n \to A_{2n}$ induced by the morphism $\id_n \oplus \iota_n: (n\oplus 0) \to (n \oplus n)$. 
The subgroup $(1 \times G_n) \subseteq G_{2n}$ acts trivially on $\rho(A_n) = A_{2n}$ by functoriality. Now, the category $\G$ is braided, and so we have natural maps from the braid group on $2n$ strands to $G_{2n}$. Then the subgroup $(1 \times G_n)\subseteq G_{2n}$ is conjugate to $(G_n \times 1)$ via the image of a braid that wraps the first $n$ strands past the second $n$ strands. It follows that $(G_n \times 1)$ also acts trivially on $A_{2n}$. Since $A_n \cong \mathrm{Res}^{G_{2n}}_{G_n \times 1} A_{2n}$ as $G_n$--representations, we conclude that $A_n$ has trivial $G_n$--action as claimed.

It follows that $G_n$ also acts trivially on $B_n \cong A_n$ for $n>d$, and so we obtain another short exact sequence of $U\G$--modules
\[ 0 \longrightarrow B \longrightarrow M(B_{d+1}) \longrightarrow D \longrightarrow 0, \]
 where $B_{d+1}$ denotes the $\G$--module that is $B_{d+1}$ in degree $0$ and zero in all degrees $>0$, and $D$ is defined as the cokernel of $B \m M(B_{d+1})$. Note that $D_n = 0$ for all $n>d$.  Condition \con{H3($k$,$a$)} implies that $\widetilde H_i(M(B_{d+1}))_n =0$ for all $n>ki +a$. The first paragraph implies that $\widetilde H_i(C)_n = \widetilde H_i(D)_n = 0$ for all $n>i+1+d$. Observe that $\widetilde C_i$ is functorial with respect to maps of $U\G$--modules, and it is an exact functor. Thus short exact sequences of $U\G$-modules give short exact sequences of central stability chains and hence a long exact sequences of central stability homology groups.   Using the long exact sequences on homology associated to these two short exact sequences, we conclude that $\widetilde H_i(A)_n = 0 $ for all $n>\max(d+i+2,ki+a)$.
\end{proof}

%\begin{corollary}\label{cor:findeg homology}
%Let $a,b,k,\ell\in \N$. Let $U\G$ be a stability category satisfying \con{H3($k$,$a$)} and \con{H4($\ell$,$b$)} with $b\ge\max(k,a)$. If $A$ is a $U\G$--module of polynomial degree $\le r$ with $r \ge 0$, then \[\widetilde H_i(A)_n = 0 \qquad \text{ whenever \quad $i\ge -1$ \quad and \quad  $n> \ell^{i+1}(r-1) + (\ell^i + \ell^{i-1} + \dots + \ell + 1)b+1$.}\] 
%\end{corollary}

By specializing \autoref{thm:findeg homology} to homological degrees $i=-1,0$ and invoking \autoref{thm:res of finite type}, we obtain the following consequences for functors of finite polynomial degree. 

%\begin{corollary} \label{cor:CSD&PresentationDeg}
%Let $a,b,k,\ell\in \N$. Let $U\G$ be a stability category satisfying \con{H3($k$,$a$)} and \con{H4($\ell$,$b$)} with $b\ge\max(k,a)$. Suppose $A$ is a $U\G$--module of polynomial degree $\le r$ in ranks $>d$ for some $r \ge 0$ and $d \geq -1$. 
%\begin{enumerate} 
%\item \label{CSDforPoly} The $U\G$--module $A$ has central stability degree $\le  \max( a, d+2, \ell (d+r) + b +1)$. 
%\item \label{PresentationDegForPoly} The $U\G$--module $A$ is generated in degrees $g \le  \max(  a-k, d+1, d+r +1)$ and $A$ is related in degrees $\le  \max( a, d+2, \ell (d+r) + b +1, g+b)$.
%\end{enumerate}
%\end{corollary}

\begin{corollary} \label{cor:CSD&PresentationDeg}
Let $a,b,k,\ell \geq 1$. Let $\G$ be a braided stability groupoid, and assume its stability category $U\G$ satisfies \con{H3($k$,$a$)} and \con{H4($\ell$,$b$)} with $b\ge\max(k,a)$.  Suppose $A$ is a $U\G$--module of polynomial degree $\le 0$ in ranks $>d$ for some $d \geq -1$. 
\begin{enumerate} 
\item \label{CSDforPoly0} Then  $A$ has central stability degree $\le  \max( a, d+2)$. 
\item \label{PresentationDegForPoly0} Then  $A$ is generated in degrees $ \le  \max(  a-k, d+1)$ and related in degrees $\le  \max( a-k+b, d+1+b)$.
\end{enumerate}
Suppose $A$ is a $U\G$--module of polynomial degree $\le r$ in ranks $>d$ for some $r \ge 1$ and $d \geq -1$. 
\begin{enumerate} 
\item \label{CSDforPoly} Then  $A$ has central stability degree $\le  \ell (d+r) + b +1$. 
\item \label{PresentationDegForPoly} Then  $A$ is generated in degrees $ \le    d+r +1$ and $A$ is related in degrees $\le   \ell (d+r) + b +1$.
\end{enumerate}
\end{corollary} 

%Moreover, $A$ admits a presentation \[ P_1 \longrightarrow P_0 \longrightarrow A\] such that $P_0$ is freely generated in ranks $\le (d+r+1)$ and $P_1$ is freely generated in ranks $\le ( \ell (d+r) + b +1)$.

%\begin{remark} It may be reasonable to express the result \autoref{cor:CSD&PresentationDeg} \ref{PresentationDegForPoly} by saying that ``$A$ is related in ranks $\le ( \ell (d+r) + b +1)$''. We note, however, that this is a distinct and stronger notion of relation degree for $\FI$--modules than is used in Church--Ellenberg \cite[Definition 4.1]{CE}, since Church and Ellenberg allow resolutions by any $\FI$--modules of the  form given in \autoref{DefnM(W)}, whereas we insist that $P_0$ and $P_1$ be direct sums of represented $\FI$--modules. 
%\end{remark} 

\subsection{$\FI$--modules of finite polynomial degree} \label{SectionFinPolyDeg}

It is easy to see that $U\mathfrak S$ is a subcategory of every symmetric stability category, in particular of $U\GL(R)$, $U\Sp(R)$, and $U\GL^H(R)$ for any $H$ containing $-1$. Note that the definition of polynomial degree never mentions the groups $G_n$ and only involves vanishing of kernels and cokernels of shift functors. In particular, it follows from the definition of polynomial degree that the polynomial degree of a $U\G$--module coincides with the polynomial degree of the underlying $\FI$--module. 

 We now compute the polynomial degree of  free $\FI$--modules.

%We now recall the definition of a collection of $\FI$--modules first introduced by Church--Ellenberg--Farb \cite[Definition 2.2.2]{CEF} and then proceed to compute their polynomial degree.  

\begin{proposition}\label{polyFI}
Let $W$ be a $\Z[\mathfrak S_m]$--module, then $M(W)$ has polynomial degree $\le m$.
\end{proposition}

\begin{proof}
We prove the assertion by induction over $m$. If $m=0$, then
\[ M(W) = M(0)\otimes_{\Z} W\]
and the map $M(W) \to SM(W)$ is an isomorphism. This implies that $M(W)$ has polynomial degree $\le 0$. Now suppose $m>0$. Church--Ellenberg \cite[Lemma 4.4]{CE} showed that
\[ \coker M(W) = M(\Res^{\mathfrak S_m}_{\mathfrak S_{m-1}} W).\]
Thus $\coker M(W)$ has polynomial degree $\le m-1$ by induction. The proof of \cite[Lemma 4.4]{CE} also shows that $\ker M(W) = 0$. This completes the induction.
\end{proof}

\begin{proposition} \label{PropH4FI}
$U\mathfrak S$ satisfies \con{H4($1$,$1$)}.
\end{proposition}

\begin{proof}
The isomorphisms
\[ \coker M(m) \cong \coker M\left(\Z[\mathfrak S_m]\right) \cong M\left(\Res^{\mathfrak S_m}_{\mathfrak S_{m-1}}\Z[\mathfrak S_m]\right) \cong M(m-1)^{\oplus m}\]
show that we can apply \autoref{thm:res of finite type}, where $W^0$ is of degree $m-1$ and $W^i=0$ for $i>0$. Invoking \autoref{H3} and choosing $d_i = m-1+i$ implies that 
\[\widetilde H_i(\coker(M(m)))_n = 0 \quad \text{for $n>d_{i+1} = m+i > 1\cdot(m+i) + 1$}.\qedhere\]
\end{proof}

\begin{theorem}\label{FIpolypres}
An $\FI$--module $A$ is presented in finite degree if and only if it has finite polynomial degree.
Specifically, if $A$ is generated in degrees $\le r$ and presented in degrees $\le d$, then $A$ has \[ \text{polynomial degree $\le r$} \quad \text{ in ranks $> d+\min(r,d) -1$.}\] If $A$ has polynomial degree $\le r$ in ranks $>d$  for some $d \geq -1$, it is 
\begin{align*} \text{generated in degrees } & \le d+r+1, \text{ and } \\   \text{presented in degrees } &\le  d+r+2. \end{align*} 
\end{theorem}

\begin{proof}
The first direction is given by Randal-Williams--Wahl \cite[Example and Proposition 4.18]{RWWarxiv}. Assume $A$ has polynomial degree $\le r$ in ranks $>d$. Then by \autoref{H3} and \autoref{PropH4FI} we may apply \autoref{cor:CSD&PresentationDeg} with the values $a=k=b=\ell=1$. 
\end{proof}

\subsection{$\VIC(R)$--, $\VIC^H(R)$--, $\SI(R)$--modules of finite polynomial degree} \label{SectionMainLemma}

In this subsection, we use the connectivity results of \autoref{sechighcon} to prove a vanishing result for the central stability homology of $\VIC^H(R)$--modules and $\SI(R)$--modules with finite polynomial degree for a $R$ a PID.  More specifically, we will show that the connectivity results of \autoref{sechighcon} establish the hypotheses of \autoref{thm:findeg homology}. We use the following criterion for \con{H4($\ell$,$b$)}. 

\begin{proposition}[Patzt {\cite[Proposition 7.10]{Pa2}}] \label{linkH4} Assume \con{H3($k$,$a$)} and let $b\ge \max(k,a)$. Then the condition
\con{H4($\ell$,$b$)} holds if for every $(m-1)$--simplex $\tau \in (K_{m-1}\G)_{n+1}$, the intersection
\[ \Lk^{(K_\bullet\G)_{n+1}}_\bullet( \tau) \cap ( K_\bullet \G)_{n} \subset (K_\bullet \G)_{n+1}\]
$\displaystyle  \left\lfloor \left(\frac{n-b}\ell \right) -m-1 \right\rfloor$--connected, where is $( K_\bullet \G)_{n} \inject  (K_\bullet \G)_{n+1}$ is induced by $0\oplus n \to 1 \oplus n$. 
\end{proposition}

\begin{proof}
The statement is a reformulation of \cite[Proposition 7.10]{Pa2}. As noted in the paragraph above \cite[Proposition 7.10]{Pa2}, $\widetilde H_*(\coker M(m))_n$ is a quotient of the associated homology of the augmented semisimplicial set $(K_\bullet T\Hom(m,-))_n$. Therefore when $(K_\bullet T\Hom(m,-))_n$ is $\left\lfloor \left( \frac{n-b}{\ell} \right) -m - 1\right\rfloor$--connected, the central stability homology $\widetilde H_i(\coker M(m))_n$ vanishes for all $n> \ell(i+m) +b$. \cite[Proposition 7.10]{Pa2} states that $(K_\bullet T\Hom(m,-))_n$ is the disjoint union of the augmented semisimplicial sets
\[  \Lk^{(K_\bullet\G)_{n+1}}_\bullet( \tau) \cap ( K_\bullet \G)_{n}\]
for $\tau \in \Hom^{U\G}(m,n+1) = (K_{m-1}\G)_{n+1}$. These augmented semisimplicial sets have singleton sets in degree $-1$. 
\end{proof}

\begin{proposition}\label{H4VIC} \label{VICH-H4}
Let $R$ be a PID of  stable rank $s$ and let $H\le R^\times$  contain $-1$. Then $U\GL^H(R)$ satisfy \con{H4($2$,$s+1$)}. 
\end{proposition}

\begin{proof} 
By \autoref{H3}, we can use \autoref{linkH4} to prove this proposition.

%$\Z[\PBC_\bullet(R^n)] \cong \widetilde C_\bullet(M(0))_n$. The second author \cite{Pa2} gave the following criterion for \con{H4($\ell$,$b$)} to hold: for every $m$--simplex $\tau=(f,C)\in \PBC_m(R^{n+1})$ the intersection
%\[ \Lk_\bullet \tau \cap \PBC_\bullet(R^n)\]
%is $\left( \frac{n-b}\ell -m-2 \right)$--connected. 

Let us first check this for $U\GL(R)$.
By \autoref{linkH4}, we must check that for every $(m-1)$--simplex $\tau \in (K_{m-1} \GL(R))_{n+1}$, the intersection
\[ \Lk^{(K_\bullet\GL(R))_{n+1}}_\bullet( \tau) \cap (K_\bullet \GL(R))_{n}\]
is $\displaystyle \left\lfloor \left(\frac{n-s-1}2\right)-m-1\right\rfloor$--connected. Recall that $\PBC_\bullet(R^n) \cong (K_\bullet \GL(R))_n$.  \autoref{linksPBC} applied to
\[  \Lk^{(K_\bullet\GL(R))_{n+1}}_\bullet( \tau) = \PBC_\bullet(R^{n+1},C,\im(f))\quad\text{and}\quad (K_\bullet\GL(R))_n = \PBC_\bullet(R^{n+1}, R^n,R ),\]where $(f,C)=\tau\in (K_{m-1}\GL(R))_{n+1} = \Hom^{U\GL(R)}(m,n+1)$ and $R \oplus R^n = R^{n+1}$, shows that the intersection appearing in \autoref{linkH4} is isomorphic to
\[ \PBC_\bullet(R^{n+1}, C \cap R^n, \sat(\im f + R) ).\]
Since 
\[ \rk (C\cap R^n) \ge n-m\]
and 
\[ \rk \sat(\im f+R) \le m+1,\]
 \autoref{PBCsemiconn} implies that the semisimplicial set is $\displaystyle \left\lfloor \frac{(n-m)-(m+1) - s-2}2 \right\rfloor$--connected. But
\[ \frac{(n-m)-(m+1) - s-2}2 = \frac{n-s-1}{2} - m -1,\]
and the result follows.   
%
%The semisimplicial $\VIC$--module $\widetilde C_i(S^r M(d))(V)$ is isomorphic to the simplicial chains on an augmented semisimplicial set whose set of $p$--simplices is given by \[ \bigsqcup_{(g,C) \in \Hom^{\VIC}(R^{p+1},V)} \Hom^{\VIC}(R^d, R^r\oplus C) \] and whose facemaps are defined in a manner similar to those of $\widetilde C_i(S^r M(d))(V)$. In \cite{Pa2} this semisimplicial set is proved to coincide with
%\[ \coprod_{(f,D)=\tau \in \PBC_{d-1}(R^r\oplus V)} \Lk\tau \cap \PBC_\bullet(V) \cong \PBC_\bullet^\alpha\]
%with $\alpha = {(R^r\oplus V, V,R^r,D,\im f)}$.
% We have shown that $\PBC_\bullet^\alpha$ is $\frac{n-2d-r-s-2}{2}$--connected in  \autoref{PBCsemiconn} because
% \[ \rk V + \rk D - \rk R^r\oplus V = n-d\]
% and 
% \[ \rk R^r\oplus V = r+n.\] Thus, the relevant homology groups vanish.
%%W=\im f, Y=F, X=V, C=Z 
%\end{proof}
%
%
%
%\begin{proposition}
%Let $R$ be a PID of  stable rank $s$, and let $H\le R^\times$. Then $U\GL^H(R)$ satisfies \con{H4($2$,$s+1$)}. 
%\end{proposition}
%
%\begin{proof} 

We will now check the assertion for $U\GL^H(R)$.
We must check that for every $m$--simplex $\tau \in (K_{m-1} \GL^H(R))_{n+1}$, the intersection
\[ \Lk^{(K_\bullet\GL^H(R))_{n+1}}_\bullet( \tau) \cap (K_\bullet \GL^H(R))_{n}\]
is $\displaystyle \left\lfloor\left( \frac{n-s-1}2\right)-m-1\right\rfloor$--connected. Because
\[ (K_p \GL^H(R))_{n+1} = (K_p \GL(R))_{n+1} \qquad \qquad \text{for $p \le n-1$, } \]
\[  \Lk^{(K_\bullet\GL^H(R))_{n+1}}_\bullet( \tau) \cap (K_\bullet \GL^H(R))_{n} \inject  \Lk^{(K_\bullet\GL(R))_{n+1}}_\bullet( \tau) \cap (K_\bullet \GL(R))_{n}\]
induces a bijection of the set of $p$--simplices for $m+p+1 < n+1$. Thus it induces an $(n-m-1)$--connected map on geometric realizations and so the assertion follows from the connectivity of the codomain which was proved in \autoref{H4VIC}.
\end{proof}

Because $U\GL^H(R)$ satisfies \con{H3($2$,$s+1$)} and \con{H4($2$,$s+1$)}, \autoref{thm:findeg homology} implies the following corollary.

\begin{corollary} \label{corpolyVIC}
Assume $R$ is a PID of  stable rank $s$ and let $H\le R^\times$ contain $-1$.
If $A$ is a $\VIC^H(R)$--module of polynomial degree $\le 0$ in ranks $>d$ for some $d \geq -1$, then \[\widetilde H_i(A)_n = 0 \qquad \text{ for all \quad $i\ge -1$ \quad and \quad  $n> \max(d+i+2, 2i+s+1)$.} \]
If $A$ is a $\VIC^H(R)$--module of polynomial degree $\le r$ in ranks $>d$ for some $r \ge 1$ and $d \geq -1$, then \[\widetilde H_i(A)_n = 0 \qquad \text{ for all \quad $i\ge -1$ \quad and \quad  $n> 2^{i+1}(d+r +s+1) -s$.} \]
In particular, $A$ has central stability degree $\le \max(d+2, s+1, 2d+2r+s+2)$.
\end{corollary}

\begin{corollary} \label{Cor:FIvsVICH}
Let $R$ be a PID of stable rank $s$, and let $H\le R^\times$ contain $-1$. Let $A$ be a $\VIC^H(R)$--module such that the underlying $\FI$--module is generated in degrees $\le g$ and related in degrees $\le r$.  Then, as a $\VIC^H(R)$--module, $A$ is  generated in degrees $\le g$ and presented in degrees
$$ \begin{array}{ll} \le \max(2s, r+\min(g,r)+ s+1) & \text{ if $g=0,$ and} \\ \le 2r+2g+2\min(g,r)+s & \text{ if $g>0$.} \end{array}$$ 
In particular, as a $\VIC^H(R)$--module, $A$ has central stability degree $$ \begin{array}{ll} \le \max(s+1, r+1) & \text{ if $g=0,$ and} \\ \le 2r+2g+2\min(g,r)+s & \text{ if $g>0$.} \end{array}$$ 
\end{corollary}

\begin{proof} Since $\FI \subseteq \VIC^H(R)$, generation in degree $\le g$ over $\FI$ implies generation in degree $\le g$ over the larger category $\VIC^H(R)$.  By \autoref{FIpolypres}, the sequence $A$ has polynomial degree $\le g$ in ranks $> r+\min(g,r)-1$, viewed either as a module over $\FI$ or over $\VIC^H(R)$. Because $\VIC^H(R)$ satisfies conditions  \con{H3($2$,$s+1$)} and  \con{H4($2$,$s+1$)}, by \autoref{H3} and \autoref{VICH-H4}, respectively, we can use \autoref{cor:CSD&PresentationDeg} to conclude the result. 
\end{proof}

\begin{proposition} \label{USpRH4}
Let $R$ be a PID of unitary stable rank $s$. Then $U\Sp(R)$ satisfies \con{H4($2$, $s+2$)}. 
\end{proposition}

\begin{proof}
We can apply \autoref{linkH4} because of \autoref{H3}. We need to check the connectivity of the intersection
\[ \Lk_\bullet \tau \cap (K_\bullet \Sp(R))_n\]
for every $(m-1)$--simplex $\tau\in (K_{m-1} \Sp(R))_{n+1} \cong \SPB_{m-1}(R^{2n+2})$. Let $U$ be the symplectic submodule of $R^{2n+2}$ generated by $\tau$. 
Then a simplex $\sigma$ is contained in $\Lk_\bullet \tau \cap (K_\bullet \Sp(R))_n$ exactly when $\sigma$ is an ordered symplectic partial basis for the $R$--module
$$ U^{\perp} \cap R^{2n} = U^{\perp}  \cap (R^2)^{\perp} = (U + R^2)^{\perp}. $$ Note that here we are viewing $R^2$ as the first two coordinates in $R^{2n+2}$ and $R^{2n}$ as the last $2n$ coordinates. We have that  $\Lk_\bullet \tau \cap (K_\bullet \Sp(R))_n = \SPB_{\bullet}(W^{\perp})$, where $W = U + R^2$. 

We can decompose $W$ as $W = U \oplus U'$, where $U'$ is a submodule of rank at most $2$ contained in $U^{\perp}$. By Mirzaii-van der Kallen \cite[Lemma 6.6]{MvdK},  we can embed $U'$ in a symplectic submodule $H$ of $U^{\perp}$ of rank $4$ as long as $\frac12 \rk(U^{\perp})=(n+1-m)$ is at least $(2+s)$. (If $(n+1-m) < (2+s)$ then our conclusion will be vacuous because then $ \left( \frac{ 1 }{2}(n-(s+2)) -m -1 \right) < -1$.) So $U\oplus H$ is a rank-$(2m+4)$ symplectic submodule of $R^{2n+2}$ containing $W$, hence $(U \oplus H)^{\perp}$ is a symplectic submodule of rank $(2n+2)-(2m+4)$ contained in $W^{\perp}$. This implies that $W^{\perp}$ has Witt index
$$g(W^{\perp}) \geq n-m-1.$$ Friedrich \cite[Theorem 3.4]{Friedrich} then implies that $\SPB_{\bullet}(W^{\perp})$ is 
$ \displaystyle \left\lfloor \frac12 \left( (n-m-1)-s-3 \right) \right\rfloor$--connected. Since  $$ \left(  \frac{ (n-m-1)-s-3}{2} \right) \; \;   \geq  \; \; \left( \frac{ n-(s+2) }{2} -m -1 \right),  $$
\autoref{linkH4} implies that  $U\Sp(R)$ satisfies \con{H4($2$, $s+2$)}, as claimed. 
\end{proof}

Combining this result with \autoref{H3} and \autoref{thm:findeg homology} gives the following Corollary. 

\begin{corollary} \label{corpolySI}
Assume $R$ is a PID of (symplectic) unitary stable rank $s$. If $A$ is a $\SI(R)$--module of polynomial degree $\le 0$ in ranks $>d$ for some $d \geq -1$, then \[\widetilde H_i(A)_n = 0 \qquad \text{ for all \quad $i\ge -1$ \quad and \quad  $n> \max(d+i+2, 2i+s+2)$.} \]
If $A$ is a $\SI(R)$--module of polynomial degree $\le r$ in ranks $>d$ for some $r \ge 1$ and $d \geq -1$, then \[\widetilde H_i(A)_n = 0 \qquad \text{ for all \quad $i\ge -1$ \quad and \quad  $n> 2^{i+1}(d+r+s+2) -s -1$.} \]

\end{corollary}

\begin{remark}
All of the results in this section and the previous section apply equally well to orthogonal groups. We chose not to include these results because we do not know of any applications. 
\end{remark}

\subsection{A spectral sequence} \label{SectionSpecSequence}

In this subsection, we summarize results about a spectral sequence introduced by  Putman and Sam. For notation that matches our paper we merge \cite[Proposition 8.2, Lemma 8.3]{Pa2} into the next proposition.

%\begin{definition}\label{Def:stability groupoid}
%We call a strict monoidal groupoid $\mathcal G$ a \emph{stability groupoid} if
%\begin{enumerate}
%\item its underlying monoid is the nonnegative numbers with addtition,
%\item there are only automorphisms,
%\item the automorphism group of $0$ is the trivial group,
%\item it is injective (see below).
%\end{enumerate}
%
%
%Let us abbreviate $G_n$ for $\Aut_{\mathcal G}(n)$. From the monoidal structure we get group homomorphisms $G_m\times G_{n-m} \to G_n$. We call $\mathcal G$ \emph{injective} if all these maps are injective. Note that we get embeddings of $G_m$ into $G_n$ if $m\le n$. Let us distinguish the \emph{left} and the \emph{right} embedding which are given by
%\[ G_m \inject G_m\times G_{n-m} \inject G_n \quad\text{and}\quad G_m \inject G_{n-m}\times G_m \inject G_n,\]
%respectively.
%\end{definition}

\begin{proposition}\label{prop:stabSES}
Let $\mathcal G, \mathcal Q$ be stability groupoids and $\mathcal G \to \mathcal Q$ a homomorphism such that $ G_n \to Q_n$ is surjective for every $n\in \mathbb N_0$. Then there is a unique stability groupoid $\mathcal N$ and a homomorphism $ \mathcal N \to \mathcal G$ such that
\[ N_n \cong \ker( G_n \longrightarrow Q_n), \]
where $N_n = \mathcal N(n)$.

Furthermore, for every $i \ge0$ there is a $U{\mathcal Q}$--module that we denote by $H_i(\mathcal N)$ with \[ H_i(\mathcal N)_n \cong H_i (N_n).\]
\end{proposition}

The following spectral sequence is central to our proof.

\begin{proposition}[Putman--Sam {\cite[Theorem 5.9]{PutmanSamLinearGroups}}, reformulated by Patzt  {\cite[Corollary 8.5]{Pa2}}]\label{prop:spectralsequence}
Let $\mathcal G, \mathcal Q$ be stability groupoids and $\mathcal G \to \mathcal Q$ a homomorphism such that $ G_n \to Q_n$ is surjective for every $n\in \mathbb N_0$ and $\mathcal N$ as in \autoref{prop:stabSES}. Assume that $\mathcal G$ and $\mathcal Q$ are braided and that $\mathcal G \m \mathcal Q$ is a map of braided monoidal groupoids. 
For each $n$,  there is a spectral sequence with
\[ E^2_{p,q} \cong \widetilde H^{\mathcal Q}_p( H_q(\mathcal N))_n .\]
Notably, $E^2_{p,q}$ vanishes for $p < -1$ or $q<0$. If $U\G$ satisfies  \con{H3($k$,$a$)} then the spectral sequence converges to zero for $ \displaystyle p+q \le \left(  \frac{n-a-1}{k} \right).$ 
\end{proposition}

%\begin{proof} This existence of these spectral sequences follows from work of Putman--Sam and the second author, and we need only verify the ranges for which they converge to zero. For a given $n$, the spectral sequence $E^r_{p,q}$ can be realized as one of the two spectral sequences arising from the double complex $E_* G_n \otimes_{RN_n} \widetilde C_*^{\mathcal G} (A)_n$. The second spectral sequence $\tilde E^r_{p,q}$ assocated to this double complex satisfies $$\tilde E^1_{p,q} \cong E_p G_n \otimes_{RN_n} \widetilde H_q^{\mathcal G} (A)_n;$$ see Patzt \cite[Proposition 8.4]{Pa2}. Thus to ensure the group $E^{\infty}_{p,q}$ is zero, it suffices to show that the diagonal $\tilde E^1_{0,p+q}, \, \tilde E^1_{1,p+q-1},  \, \ldots,  \,   \tilde E^1_{p+q,0},  \tilde E^1_{p+q+1,-1}$ vanishes. Condition  \con{H3($k$,$a$)} implies a range of values of $n$ for which these groups vanish, with the most restrictive bound imposed by the vanishing of $\tilde E^1_{0,p+q}$ when $n > k(p+q)+a$.  \end{proof} 

%%%%%%%%%%%%% SECTION %%%%%%%%%%%%%%%%%%%%

\section{Applications} 

\label{secap}
In this section, we prove our central stability results for the second homology groups of Torelli subgroups and for the homology of congruence subgroups.

\subsection{$H_2(\IA_n)$} 

We will recall a computation of $H_1(\IA_n)$. Then we will use a spectral sequence argument and our results about polynomial $\VIC(R)$--modules to prove \autoref{IAn}, central stability for $H_2(\IA_n)$.

Recall that $\AutF$ is the stability groupoid given by $\AutF_n = \Aut(F_n)$. Abelianization induces a monoidal functor \[ \AutF \longrightarrow \GL\] that is surjective in every degree. By \autoref{prop:stabSES}, the kernels form a stability groupoid which we denote by $\IA$. 

Reinterpreting the work of Andreadakis \cite[Section 6]{Andr}, Cohen--Pakianathan (unpublished), Farb (unpublished) and Kawazumi \cite[Theorem 6.1]{Kawazumi} in the language of $\FI$--modules, Church, Ellenberg, and Farb \cite[Equation (25) and the proof of Theorem 7.2.3]{CEF} proved the following. See also Day--Putman \cite[Page 5]{DP}.

\begin{theorem}[Andreadakis, Farb, Kawazumi, Cohen--Pakianathan, Church--Ellenberg--Farb] Use the following notation for bases for the $\FI$--module $M(1)$: 
$$M(1)_n \cong \Z[ e_1, e_2, \ldots, e_n ].$$ Then
\begin{align*} H_1(\IA) & \cong \left( {\bigwedge}^2 M(1) \right ) \otimes M(1) \\
&  \cong M\Big( \Z [ e_1 \wedge e_2 \otimes e_1, \;  e_1 \wedge e_2 \otimes e_2 ] \Big) \oplus M\Big( \Z [ e_1 \wedge e_2 \otimes e_3, \;  e_1 \wedge e_3 \otimes e_2, \;  e_2 \wedge e_3 \otimes e_1 ] \Big) 
\end{align*}
 as an $\FI$--module. 
 %Rationally, \[ H_1(\IA; \Q) \cong M(2) \oplus M\left(\tiny\yng(2,1)\right) \oplus M\left(\tiny\yng(1,1,1)\right).\]
 \label{H1IAn}  \end{theorem}
 
 %%%%%%%%%%%%% START PURPLE
 We are now ready to prove \autoref{IAn}  and the following variant \hyperlink{thmA'}{Theorem A'}.  
 \begin{ThmA'} \hypertarget{thmA'}{}
The $\VIC(\Z)$--module $H_2(\IA)$ is generated in degrees $\leq 18$ and presented in degrees $\leq 38$. 
 \end{ThmA'}

%We now prove \autoref{IAn}, central stability for $H_2(\IA)$.

\begin{proof}[Proof of  \autoref{IAn} and \hyperlink{thmA'}{Theorem A'}]
$\AutF$ is symmetric monoidal and thus so is $U{\AutF}$. This follows from  \autoref{prop:symmetric monoidal}. We can therefore apply the spectral sequences of  \autoref{prop:spectralsequence} with $\mathcal N = \IA$, $\mathcal G = \AutF$, $\mathcal Q= \GL(\Z)$.
%and $A$ the trivial $U\AutF$--module $M(0)$.
\autoref{H3} says that $U\AutF$ satisfies \con{H3($2$,$2$)}, hence these spectral sequences converge to zero for $p+q\le \frac{n-3}2$. 

We will show that 
\[E^2_{-1,2} \cong 0 \quad\text{for $n> 18$}\qquad\text{and}\qquad E^2_{0,2} \cong \widetilde H^{\GL(\Z)}_0( H_2(\IA))_n = 0 \quad \text{for $n > 38$}\]
by showing that there are no nontrivial differentials to or from the groups $E^2_{-1,2}$ and $E^2_{0,2}$ of the spectral sequence in this range, as in  \autoref{E2PageVanishing}. 
%%%%%%
%%%%%%
%%%%%%
\begin{figure}[h!]    \centering \begin{tikzpicture} \footnotesize 
  \matrix (m) [matrix of math nodes,
    nodes in empty cells,nodes={minimum width=3ex,
    minimum height=5ex,outer sep=2pt},
    column sep=6ex,row sep=3ex]{ 
              %  &      &     &     &\strut& & & &  \\    
 3    &   \bigstar &  \bigstar &  \bigstar &\bigstar & \bigstar  &  \\  
2    &  \widetilde H^{\GL(\Z)}_{-1}( H_2(\IA))_n & \widetilde H^{\GL(\Z)}_0( H_2(\IA))_n & \widetilde H^{\GL(\Z)}_1( H_2(\IA))_n  &  \widetilde H^{\GL(\Z)}_2( H_2(\IA))_n  &  \\           
1     & 0 & 0 &0  & 0 &   \bigstar  \\           
 0     &  0  & 0   & 0   &0 &0  &  \\       
 \quad\strut &   -1  &  0  &  1  & 2  &3  &   \\}; 
 
 \draw[-stealth, blue] (m-3-4.west) -- (m-2-2.east) node [midway,above] {$d_2$};
 \draw[-stealth, red] (m-4-5.west) -- (m-2-2.east) node [midway,below] {$d_3$};
 
 \draw[-stealth, blue] (m-3-5.west) -- (m-2-3.east) node [midway,above] {$d_2$};
 \draw[-stealth, red] (m-4-6.west) -- (m-2-3.east) node [midway,below] {$d_3$};

\draw[thick] (m-1-1.east) -- (m-5-1.east) ;
\draw[thick] (m-5-1.north) -- (m-5-6.north) ;

 %\draw [line width=0.05mm, lightgray] (m-1-3) -- (m-5-3) ;
 %\draw [line width=0.05mm, lightgray] (m-4-1) -- (m-4-6) ;
\end{tikzpicture}
\caption{Page $E^2_{p,q}$ for $n \ge 39$.} \label{E2PageVanishing}
\end{figure}
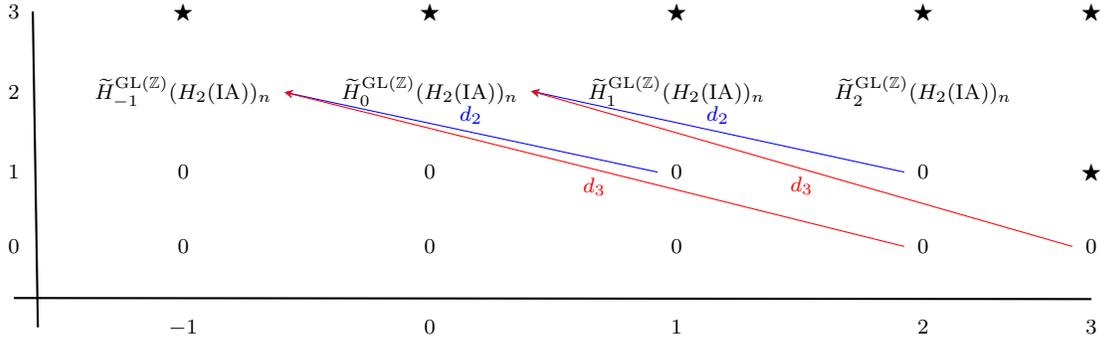  
%%%%%%
%%%%%%
%%%%%%
By \autoref{H3}, the central stability homology of the $\VIC(\Z)$--module $H_0(\IA) \cong M(0)$ is zero in degree $q$ for $n > 2q+3$. This implies the vanishing of the groups $E^2_{2,0} \cong \widetilde H^{\GL(\Z)}_{2}( H_0(\IA))_n$ and $E^2_{3,0} \cong \widetilde H^{\GL(\Z)}_{3}( H_0(\IA))_n$ for $n > 9$.  \autoref{H1IAn} and \autoref{polyFI} imply that $H_1(\IA)$ has polynomial degree $\le 3$.
Thus, by  \autoref{corpolyVIC}, its central stability homology as a $\VIC(\Z)$--module $\widetilde H^{\GL(\Z)}_q(H_1(\IA))_n$ vanishes for $n > 2^{q+1}(5)-2$. 
It follows that the groups $E^2_{1, 1} \cong \widetilde H^{\GL(\Z)}_{1}( H_1(\IA))_n$  vanish for $n> 2^{1+1}(5)-2=18$ and that the groups $E^2_{2,1} \cong \widetilde H^{\GL(\Z)}_{2}( H_1(\IA))_n$ vanish for $n> 2^{2+1}(5)-1=38$. Thus $E^2_{-1,2} \cong 0$ for $n>18$ and $E^2_{0,2} \cong 0$ for $n> 38$. 

Since $E^2_{-1,2} \cong \widetilde H^{\GL(\Z)}_{-1}( H_2(\IA))_n$, \autoref{lem:generation} implies that $H_2(\IA)$ is generated in degrees $\le 18$. Because $E^2_{0,2} \cong \widetilde H^{\GL(\Z)}_0( H_2(\IA))_n$, we deduce that $H_2(\IA)$ has  centrally stability degree $\le 38$. Using \autoref{thm:centralpresent}, this implies that $H_2(\IA)$ is presented in degrees $\le 38$.
\end{proof}

\subsection{$H_2(\mathcal I_g)$}

Using a similar proof strategy to that of the previous subsection, we will prove \autoref{Ig}, central stability for $H_2(\mathcal I)$.

%Its monoidal structure  and a braiding is described in \cite[Sec 5.6]{RWW}. 

Recall that $\Mod$ is the stability groupoid given by the mapping class groups of compact, connected, oriented surfaces with one boundary component. Its action on the first homology group of the surface (which preserves the intersection form) gives a monoidal functor
\[ \Mod \longrightarrow \Sp\]
that is surjective in every degree. Let $\mathcal I$ denote the stability groupoid given by assembling the kernels of this map using \autoref{prop:stabSES}.
%\end{definition}

Similar to \autoref{H1IAn}, we have the following description of $H_1(\mathcal I)$.

\begin{theorem}[Johnson {\cite[Theorem 3]{JohnsonAbelianization}; see Brendle--Farb \cite[Sections 1 and 2.2]{BrendleFarb}}] \label{PropJohnsonAbelianization}  As a module over $\FI \subseteq \SI(\Z)$, the first homology groups $H_1(\mathcal I)$ decompose as follows. 
Let \[H_g:= H_1(\Sigma_{g,1};\Z) \cong M(1)_g^{\oplus 2}, \qquad  H_g \cong \Z[ a_1, b_1, a_2, b_2, \ldots, a_g, b_g].\] Then for all $g \geq 3$, 
\begin{align*} H_1(\mathcal{I}_g ) \cong &  \phantom{\oplus} {\bigwedge}^3 H_g \oplus  \Big(   \mathrm{Sym}^0(H_g)  \oplus \mathrm{Sym}^1(H_g) \oplus \mathrm{Sym}^2(H_g)/H_g  \Big) \otimes_{\Z}  \Z / 2\Z 
\end{align*}
and these isomorphisms respect the $\FI$--module structure. In particular, in the range $g \geq 3$, the $\FI$--module $H_1(\mathcal{I} )$ coincides with $\FI$--module
\begin{align*}
 %\bigg(
& \phantom{\oplus} M\Big( \Z[ \; a_i \wedge b_i \wedge b_j, \; a_i\wedge a_j \wedge b_i \; | \; \{i,j\} = \{1,2\} ] \Big) \\ 
& \oplus M\Big( \Z[ \; a_i \wedge a_j \wedge a_k, \; a_i\wedge a_j \wedge b_k, \; a_i\wedge b_j \wedge b_k, \; b_i\wedge b_j \wedge b_k \; | \; \{i,j,k\} = \{1,2,3\} ] \Big) %\bigg) 
\\
& \oplus% \bigg( 
 M \left(\Z / 2\Z[\mathfrak{S}_0] \right)  \oplus M\left( \Z / 2\Z[ a_1, b_1, a_1b_1] \right) \oplus M\left( \Z / 2\Z[ a_1a_2, \; a_1b_2, \; a_2b_1,  \; b_1b_2  ] \right).% \bigg)
\end{align*}
%Rationally, \begin{align*} H_1(\mathcal{I} ; \Q) \cong    \left( M(2)^{\oplus 2} \oplus M\left(\Y{1,1,1}\right)^{\oplus 4} \oplus M\left(\Y{2,1}\right)^{\oplus 2} \right) 
%\oplus  \Big(  M(0) \oplus M(1)^{\oplus 3}  \oplus  M(\Y{2})^{\oplus 2}  \oplus M(2)  \Big)    \otimes_{\Z}  \Z/2\Z \end{align*} 
\end{theorem}

This description of the abelianization $H_1(\mathcal{I}_{g,1}; \Z)$ of the Torelli group $\mathcal{I}_{g,1}$ was first computed by Johnson \cite[Theorem 3]{JohnsonAbelianization} in a series of papers \cite{JohnsonQuadratic, JohnsonGenerators, JohnsonBP, JohnsonAbelianization} building on work of Birman--Craggs \cite{BirmanCraggs}. See Brendle--Farb \cite[Sections 1 and 2.2]{BrendleFarb}  or van den Berg \cite[Theorem 3.5.6]{vdBerg} for a more self-contained statement of the description of $H_1(\mathcal{I}_g )$. 
This isomorphism arises as follows. The \emph{Johnson homomorphism} is the map 
$$ \tau_g: \cI_g \twoheadrightarrow {\bigwedge}^3 H_g$$ 
which admits several algebraic and topological descriptions;  see for example Farb--Margalit \cite[Section 6.6]{FarbMargalitPrimer}. The maps $\tau_g$ are equivariant with respect to the action of the mapping class group, and natural with respect to the maps $\cI_k \to \cI_g$ and $H_k \to H_g$ induced by the embedding of a subsurface $\Sigma_{k,1} \hookrightarrow \Sigma_{g,1}$. 
The torsion subgroup of the abelianization
$$B_g^2 := \Big(   \mathrm{Sym}^0(H_g)  \oplus \mathrm{Sym}^1(H_g) \oplus \mathrm{Sym}^2(H_g)/H_g  \Big) \otimes_{\Z}  \Z / 2\Z $$ 
is identified with the space of \emph{Boolean polynomials} of degree at most $2$ in the elements $\overline{a_1}, \overline{b_1},  \ldots, \overline{a_g}, \overline{b_g}$. These elements represent functions $\Omega \to \Z/2\Z$ from the space $\Omega$ of mod-2 self-linking forms on $H_g$ associated to embeddings of $\Sigma_{g,1}$ into $S^3$. 
 The map
$$ \sigma_g: \cI_g \twoheadrightarrow B_g^2$$ 
is the \emph{Johnson--Birman--Craggs homomorphism} (after modding out by Boolean polynomials in homogeneous degree 3). The maps $\sigma_g$ are also equivariant with respect to the action of the mapping class group and natural with respect to embeddings $\Sigma_{k,1} \hookrightarrow \Sigma_{g,1}$. We caution, however, that the action of the mapping class group (and induced action of $\Sp_{2g}(\Z)$) on the elements  $\overline{a_1}, \overline{b_1},  \ldots, \overline{a_g}, \overline{b_g}$ does not simply correspond to its action on $H_g$; see for example Johnson \cite[Section 2]{JohnsonAbelianization} or Brendle--Farb \cite[Section 2]{BrendleFarb} for an explicit description of the space $B_g^2$ and map $\sigma_g$. The action of the symmetric group $\mathfrak S_g \subseteq \Sp_{2g}(\Z)$ on $B_g^2$  is, however, induced by its standard permutation action on $H_g$, extended diagonally to monomials in the elements $\overline{a_i}, \overline{b_i}$. The $\FI$ action on the torsion subgroup of $H_1(\mathcal{I}_g )$ for $g \geq 3$ is induced by the usual action on $H=M(1)^{\oplus 2}$. 

 %The free part,  $ \bigwedge^3 H$, is the image of the \emph{Johnson homomorphism}  \cite{JohnsonQuotient}. The torsion part can be identified with a certain space of \emph{Boolean polynomials} of degree at most $2$ with coefficients in $\Z/2\Z$; these polynomials represent part of the dual space to the space of mod 2 self-linking forms on $H$. 
%See for example Farb--Margalit \cite[Section 6.6]{FarbMargalitPrimer} and Brendle--Farb \cite[Section 2]{BrendleFarb} for a modern exposition. 

%We now prove \autoref{Ig}, central stability for $H_2(\mathcal I)$.

 %%%%%%%%%%%%%%%%%%%%%%%%%%% START PURPLE
 We can now prove \autoref{Ig} and the following variant \hyperlink{thmB'}{Theorem B'}.  
 \begin{ThmB'} \hypertarget{thmB'}{}
The $\SI(\Z)$--module $H_2(\mathcal{I})$ is generated in degrees $\leq 33$ and is presented in degrees $\leq 69$. 
 \end{ThmB'}

\begin{proof}[Proof of  \autoref{Ig} and \hyperlink{thmB'}{Theorem B'}]
We can apply  \autoref{prop:spectralsequence} with $\mathcal N = \mathcal I$, $\mathcal G = \Mod$, and $\mathcal Q= \Sp(\Z)$.

\autoref{H3} implies that $U\Mod$ satisfies \con{H3($2$,$2$)}, so the spectral sequence converges to zero for $p+q \le \frac{n-3}2$. By  \autoref{PropJohnsonAbelianization} and \autoref{polyFI}, $H_1(\mathcal I)$
 has polynomial degree $\le 3$ in ranks $>2$, thus by  \autoref{corpolySI} its central stability homology $\widetilde H_q^{\Sp(\Z)}(H_1(\mathcal I))_n$ vanishes for $n > 2^{q+1} \cdot 9 -3$. Moreover $\widetilde H_q^{\Sp(\Z)}(H_0(\mathcal I))_n \cong \widetilde H_q^{\Sp(\Z)}(M(0))_n$ vanishes for $n> 2q+4$. This implies that $E^2_{-1,2}=0$ for $n >  2^{1+1}\cdot 9-3=33$ and $E^2_{0,2}=0$ for $n> 2^{2+1}\cdot 9-3=69$. This proves that $H_2(\mathcal I)$ is generated in degrees $\leq 33$ and has central stability degree $\le 69$. The bound on its degree of presentation follows from \autoref{thm:centralpresent}. 
% \todo{check ranges}. 
\end{proof}

%%%%%%%%%%%%%%%%%%%%%%%%%%%%%%%%%%%% END PURPLE

\subsection{Congruence subgroups}

Let $R$ be a ring and $I\subseteq R$ be a two-sided ideal. Let $\U$ denote the image $R^\times$ in $R/I$. The quotient map $R \surject R/I$ induces a homomorphisms
\[ \GL_n(R) \m \GL_n^\U(R/I).\] We denote the kernel by $\GL_n(R,I)$, which is called the \emph{level $I$ congruence subgroup} of $\GL_n(R)$. These assemble to form an $\FI$--module \cite[Application 1]{CEFN}. Building the work of Putman \cite[Theorems B and C]{Pu}, Church--Ellenberg--Farb--Nagpal \cite[Theorem D]{CEFN}, Church--Ellenberg \cite[Theorem D']{CE}, Church--Miller--Nagpal--Reinhold \cite[Application B]{CMNR}, Gan and Li proved the following representation stability result for congruence subgroups viewed as $\FI$--modules.

\begin{theorem}[Gan--Li {\cite[Theorem 11]{GanLiCongruenceSubgroups}}]\label{CEcong} %Let $\mathcal{H}_i$ denote the $\FI$--module $ H_i(\GL(R,I))$. 
Let $R$ be a ring of stable rank $t$ and let $I$ be a two-sided ideal of $R$. Then for $i \geq 0$,  $H_i(\GL(R,I))$ is generated in degrees \[\le 4i+2t-1\] and presented in degree \[ \le 4i+2t+4\] when viewed as an $\FI$--module.
\end{theorem}

We have three remarks on the statement of Gan--Li {\cite[Theorem 11]{GanLiCongruenceSubgroups}}. First note that they use a different indexing convention for stable rank; their $d$ is our $(t-1)$. Secondly, they use a different definition of \emph{presentation degree} than we do here, however, the two definitions are equivalent (see for example Church--Ellenberg \cite[Proposition 4.2]{CE} for one direction; the other direction follows from \cite[Lemma 2.3]{CE} and a computation of the $\FI$--homology groups using an acyclic resolution). Third, note that although \cite[Theorem 11]{GanLiCongruenceSubgroups} assumes $I$ is proper, their theorem is also true in the case that $I=R$. Since $\GL_n(R,R) \cong \GL_n(R)$, van der Kallen's homological stability result \cite{V} implies representation stability (see \cite[Proof of Application B]{CMNR}).

We will now restrict attention to the case that the quotient map $\GL_n(R) \m \GL_n^\U(R/I)$ is surjective. This gives a surjective homomorphism of stability groupoids
\[ \GL(R) \surject \GL^\U(R/I).\]  By \autoref{prop:stabSES}, these groups assemble to form a stability groupoid which we denote by $\GL(R,I)$. \autoref{CEcong} and \autoref{Cor:FIvsVICH} together imply the following result which in turn implies \autoref{GL(R,I)}. Recall that the group $\GL_n(R,I)$ depends only on $I$ as a nonunital ring, and not on the ambient ring $R$. 

\begin{ThmC'}

% Let $R$ be a commutative ring and $I \subsetneq R$ a nonzero ideal, such that $R/I$ is a PID of stable rank $s$. Let $t$ be a bound on the stable rank of any commutative unital ring containing $I$. 

Let $I$ be a two-sided ideal in a ring $R$ of stable rank $r$ and let $t$ be the minimal stable rank of all rings containing $I$ as a two-sided ideal. Let $\U$ denote the image of $R^\times$ in $R/I$. If $R/I$ is a PID of stable rank $s$, and the natural map $ \GL_n(R) \m \GL^{\U}_n(R/I)  $ is surjective,  then the $\VIC^{\mathfrak{U}}(R/I)$--module  $H_i(\GL(R,I))$ is generated in degrees  \[\begin{array}{lr}
\le 0 & \text{ for } i=0 \\
 \le  \min( 3 +2t, \, \max(3+r,3+s) ) & \text{ for } i=1 \\
 \le 4i +2t -1 & \text{ for } i>1
\end{array} \]  and is presented in degrees  \[\begin{array}{lr}
\le s+1 & \text{ for } i=0 \\
 \le  \min( 28 +12t +s, \, \max(5+r,5+s) ) & \text{ for } i=1 \\
 \le 24i +12t +s +4 & \text{ for } i>1.
\end{array} \]

\end{ThmC'}

 %%%%%%%%%%%%%%%%%%%%%%%%%%%%%%%%%%%%%% END PURPLE 

\begin{proof}[Proof of \autoref{GL(R,I)} and \hyperlink{thmC'}{Theorem C'}]

%Let $t$ be the stable rank of the smallest ring containing $I$ as an ideal and $s$ be the stable rank of $R/I$. 

Gan and Li's \autoref{CEcong} states bounds on the generation and presentation degree of underlying $\FI$--module structure on $H_i(\GL(R,I))$, and  \autoref{Cor:FIvsVICH}  allows us to deduce bounds on the generation and relation degree for the $\VIC^{\mathfrak{U}}(R/I)$--module structure. Then (because the $\FI$--module structure on $H_i(\GL(R,I))$ does not depend on the choice of the ring $R$), the bounds for 
$i \ge 2$ and the first term in the bound on $i=1$ follow.  For $i=0$, the $\VIC^{\mathfrak{U}}(R/I)$--module $H_0(\GL(R,I))=M(0)$  is generated in degrees $0$. It is presented in degrees $\leq (s+1)$ by \autoref{H3}. For the second term in the case $i=1$, we use the spectral sequence argument of the previous two subsections. \autoref{thm:res of finite type} allows us to relate bounds on generation and presentation degree to bounds on the central stability degree of $H_i(\GL(R,I))$. 
\end{proof}

\bibliographystyle{amsalpha}
\bibliography{SecondHomologyGroups}

\end{document}